\documentclass[12pt, reqno]{amsart}
\usepackage{amsmath, amsthm, amscd, amsfonts, amssymb, graphicx, color}
\usepackage[bookmarksnumbered, colorlinks, plainpages,  bookmarks,%
backref=page,%
pagebackref]{hyperref}

\PII{\relax}
\copyrightinfo{}{\relax}

\allowdisplaybreaks[1]
\usepackage{inputenc}
\usepackage[english]{babel}
\usepackage[abbrev,backrefs,msc-links,nobysame]{amsrefs}
\BibSpecAlias{incollection}{inproceedings}
\BibSpec{article}{%
    +{}  {\PrintAuthors}                {author}
    +{,} { \emph }                            {title}
    +{.} { }                            {part}
    +{:} { }                            {subtitle}
    +{.} { \PrintContributions}         {contribution}
    +{.} { \PrintPartials}              {partial}
    +{,} { }                       {journal}
    +{}  { \textbf}                            {volume}
    +{} { \parenthesize}                 {date}
    +{, no.}  { }               {number}
    +{,} { }                             {pages}
    +{,} { }                            {status}
    +{.} { \PrintTranslation}           {translation}
    +{.} { Reprinted in \PrintReprint}  {reprint}
    +{.} { }                            {note}
    +{.} {}                             {transition}
}

\BibSpec{partial}{%
    +{}  {}                             {part}
    +{:} {  \emph}                            {subtitle}
    +{.} { \PrintContributions}         {contribution}
    +{,} { }                       {journal}
    +{}  { }                            {volume}
    +{}  { \parenthesize}               {number}
    +{:} {}                             {pages}
    +{,} { \PrintDateB}                 {date}
}

\BibSpec{book}{%
    +{}  {\PrintPrimary}                {transition}
    +{,} { \emph}                       {title}
    +{,} { }                            {part}
    +{:} { \emph}                       {subtitle}
    +{,} { }                            {series}
    +{,} { \voltext}                    {volume}
    +{.} { Edited by \PrintNameList}    {editor}
    +{.} { Translated by \PrintNameList}{translator}
    +{.} { \PrintContributions}         {contribution}
    +{,} { }                            {publisher}
    +{,} { }                            {organization}
    +{,} { }                            {address}
    +{,} { \PrintEdition}               {edition}
    +{,} { \PrintDateB}                 {date}
    +{.} { }                            {note}
    +{.} {}                             {transition}
    +{.} { \PrintTranslation}           {translation}
    +{.} { Reprinted in \PrintReprint}  {reprint}
    +{.} {}                             {transition}
}

\BibSpec{collection.article}{%
    +{}  {\PrintAuthors}                {author}
    +{,} { \emph}                            {title}
    +{.} { }                            {part}
    +{:} { }                            {subtitle}
    +{.} { \PrintContributions}         {contribution}
    +{.} { \PrintConference}            {conference}
    +{.} { \PrintBook}                  {book}
    +{.} { In \PrintEditorsA}              {editor}
    +{.} {  }                         {booktitle}
    +{,} { }                      {pages}
    +{,} { }                            {publisher}
    +{,} { }                            {organization}
    +{,} { }                            {address}
    +{,} { \PrintDateB}                 {date}
    +{.} { \PrintTranslation}           {translation}
    +{.} { Reprinted in \PrintReprint}  {reprint}
    +{.} { }                            {note}
    +{.} {}                             {transition}
}

 \BibSpec{inproceedings}{%
     +{}  {\PrintAuthors}                {author}
     +{,} {  \emph }                            {title}
     +{.} { }                            {part}
     +{:} { }                            {subtitle}
     +{.} { \PrintContributions}         {contribution}
     +{.} { \PrintConference}            {conference}
     +{.} { \PrintBook}                  {book}
     +{.} { In \PrintEditorsA}              {editor}
     +{.} { }                         {booktitle}
     +{,} { }                      {pages}
     +{,} { }                            {publisher}
     +{,} { }                            {organization}
     +{,} { }                            {address}
     +{,} { \PrintDateB}                 {date}
     +{.} { \PrintTranslation}           {translation}
     +{.} { Reprinted in \PrintReprint}  {reprint}
     +{.} { }                            {note}
     +{.} {}                             {transition}
 }

\BibSpec{conference}{%
    +{ }  {  \emph}                        {title}
    +{ }  {\PrintConferenceDetails} {transition}
}

\BibSpec{innerbook}{%
    +{,} { \emph}                       {title}
    +{.} { }                            {part}
    +{:} { \emph}                       {subtitle}
    +{.} { }                            {series}
    +{,} { \voltext}                    {volume}
    +{.} { Edited by \PrintNameList}    {editor}
    +{.} { Translated by \PrintNameList}{translator}
    +{.} { \PrintContributions}         {contribution}
    +{.} { }                            {publisher}
    +{.} { }                            {organization}
    +{,} { }                            {address}
    +{,} { \PrintEdition}               {edition}
    +{,} { \PrintDateB}                 {date}
    +{.} { }                            {note}
    +{.} {}                             {transition}
}

\BibSpec{report}{%
    +{}  {\PrintPrimary}                {transition}
    +{,} { \emph}                       {title}
    +{.} { }                            {part}
    +{:} { \emph}                       {subtitle}
    +{.} { \PrintContributions}         {contribution}
    +{.} { Technical Report }           {number}
    +{,} { }                            {series}
    +{.} { }                            {organization}
    +{,} { }                            {address}
    +{,} { \PrintDateB}                 {date}
    +{.} { \PrintTranslation}           {translation}
    +{.} { Reprinted in \PrintReprint}  {reprint}
    +{.} { }                            {note}
    +{.} {}                             {transition}
}

\BibSpec{thesis}{%
    +{}  {\PrintAuthors}                {author}
    +{,} { \emph}                       {title}
    +{:} { \emph}                       {subtitle}
    +{.} { \PrintThesisType}            {type}
    +{.} { }                            {organization}
    +{,} { }                            {address}
    +{,} { \PrintDateB}                 {date}
    +{.} { \PrintTranslation}           {translation}
    +{.} { Reprinted in \PrintReprint}  {reprint}
    +{.} { }                            {note}
    +{.} {}                             {transition}
}
\numberwithin{equation}{section}

\makeatletter
  \def\p@enumi{\thethm}

\makeatother
\providecommand{\Space}[3][]{\ensuremath{\mathbb{#2}^{#3}_{#1}{}}}
\providecommand{\FSpace}[3][]{\ensuremath{\ifx#2l \ell_{#3}^{#1}{}\else
    #2_{#3}^{#1}{}\fi}} 
\providecommand{\rmi}{\mathrm{i}}
\providecommand{\oper}[1]{\mathcal{#1}}
\providecommand{\uir}[3][0]{\ifcase #1{\rho^{#2}_{#3}}%
\or {\breve{\rho}^{#2}_{#3}}%
\or {\tilde{\rho}^{#2}_{#3}}%
\or {\hat{\rho}^{#2}_{#3}}\fi}
\providecommand{\algebra}[1]{\ensuremath{\mathfrak{#1}}}
\providecommand{\SL}[1][2]{\ensuremath{\mathrm{SL}_{#1}(\Space{R}{})}}
\providecommand{\SU}[1][1,1]{\ensuremath{\mathrm{SU}(#1)}}
\providecommand{\norm}[2][\relax]{\left\|#2\right\|\ifx#1\relax\else_{#1}\fi}
\providecommand{\modulus}[2][\relax]{\left| #2 \right|\ifx#1\relax\else_{#1}\fi}
\providecommand{\map}[1]{\mathsf{#1}}
\providecommand{\such}{\,\mid\,}
\providecommand{\myhbar}{\hbar}
\providecommand{\myh}{h}
\providecommand{\MR}[1]{MR\href{http://www.ams.org/mathscinet-getitem?mr=#1}{#1}}
\providecommand{\Zbl}[1]{Zbl\href{http://www.emis.de:80/cgi-bin/zmen/ZMATH/en/zmathf.html?first=1&maxdocs=3&type=html&an=#1&format=complete}{#1}}
\providecommand{\href}[2]{#2}
\providecommand{\amscite}[3]{\cite{#1}#2{#3}}

\providecommand{\scalar}[3][\relax]{\left\langle #2,#3 
        \right\rangle\ifx#1\relax\else_{#1}\fi}
\providecommand{\myeprint}[2]{E-print: \href{#1}{\texttt{#2}}}

\newtheorem{thm}{Theorem}[section]
\newtheorem{prop}[thm]{Proposition}
\newtheorem{lem}[thm]{Lemma}
\newtheorem{cor}[thm]{Corollary}
\theoremstyle{definition}
\newtheorem{defn}[thm]{Definition}
\newtheorem{example}[thm]{Example}
\theoremstyle{remark}
\newtheorem{rem}[thm]{Remark}

\providecommand{\cites}[1]{\cite{#1}}
\providecommand{\citelist}[1]{#1}
\makeindex
\begin{document}

\setcounter{page}{1}


\title[Relative Convolutions
 and Covariant Transform]{Calculus of Operators:\\ Covariant Transform
 and Relative Convolutions}

\author{Vladimir V. Kisil}
\thanks{On leave from the Odessa University.}

\address{School of Mathematics,
University of Leeds,
Leeds LS2\,9JT,
England}
\email{\href{mailto:kisilv@maths.leeds.ac.uk}{kisilv@maths.leeds.ac.uk}}

\begin{abstract}
  The paper outlines a covariant theory of operators related to 
  groups and homogeneous spaces. A methodical use of groups and their
  representations allows to obtain results of algebraic and analytical
  nature. The consideration is systematically illustrated by a
  representative collection of examples.
\end{abstract}


\maketitle
\tableofcontents

\section{Introduction}
\label{sec:introduction}

Calculus of operators on groups and homogeneous spaces has a long
history~\cites{FollStein82,Dynin75,Dynin76,Howe80a,Howe80b,HoweRatcliffWildberger84,Ratcliff85a,Folland94,MTaylor84,Kisil93b,Kisil94a,Kisil94e,Kisil98a,Kisil96e,Kisil94f}
and is enjoying a recently revived
interest~\cites{Street08a,RuzhanskyTurunen10a,FischerRuzhansky12a,BahouriFermanian-KammererGallagher12,Kisil12b,Wildberger05b}.
There are some missing connections between two periods and the purpose
of this presentation is to bridge the gap. 

\section{Calculus of Pseudodifferential Operators}

The theory of pseudo-differential operators%
\index{pseudo-differential operator|see{PDO}}%
\index{operator!pseudo-differential|see{PDO}} (PDO\index{PDO}) is an
important and profound area of analysis with numerous
applications~\cites{Hormander85,Shubin87,MTaylor81,HWeyl,deGosson11a}. In
the simplest one-dimensional case, a PDO \(A=a(x,D)\) is defined from
its \emph{Weyl symbol}%
\index{Weyl!symbol}%
\index{symbol!Weyl} \(a(x,\xi)\)---a function on \(\Space{R}{2}\)---by
the identity~\cite{Folland89}*{(2.3)}:
\begin{equation}
  \label{eq:PDO-Weyl-defn}
    [Au](y)=\int_{\Space{R}{}\times\Space{R}{}}
   a({\textstyle\frac{1}{2}}(y+x),\xi)\,e^{2\pi\rmi(y-x)\xi}\,u(x)\,dx\,d\xi.
\end{equation}
The alternative \emph{Kohn--Nirenberg correspondence}%
\index{Kohn--Nirenberg symbol}%
\index{symbol!Kohn--Nirenberg} between symbols and
operators~\cite{Folland89}*{\S~2.2} is provided by a similar formula:
\begin{equation}
  \label{eq:PDO-KN-defn}
  [A_{\text{KN}}u](y)=\int_{\Space{R}{}\times\Space{R}{}}
  a(y,\xi)\,e^{2\pi\rmi(y-x)\xi}\,u(x)\,dx\,d\xi.
\end{equation}

There is a natural demand to generalise PDO for other settings.
It is common to have several competing approaches for this. We briefly
outline two of them.

\subsection{Pontryagin Duality and the Fourier Transform}
\label{sec:pontry-dual-four}

Historically, the theory of PDO\index{PDO} grown out of the study of
singular integral operators%
\index{singular!integral operator|see{SIO}}%
\index{operator!singular integral|see{SIO}} (SIO\index{SIO}), which
can be viewed either as convolutions on the Euclidean group or Fourier
multipliers. In either case, this prompts a consideration of groups and
representation theory. For simplicity, we take \(G=(\Space{R}{},
+)\)---the abelian group of reals with addition. Its Pontryagin
dual---the collection of all unimodular characters \(\chi_\xi(x)=
e^{2\pi\rmi \xi x}\)---is again the abelian group \(\hat{G}\)
isomorphic to \((\Space{R}{},+)\)~\amscite{KirGvi82}*{\S~IV.2.1}. The
Fourier transform \(\oper{F}\) maps a function \(f\) from the Schwartz
space%
\index{Schwartz!space}%
\index{space!Schwartz} \(\FSpace{S}{}(G)\) to \(\hat{f}\in
\FSpace{S}{}(\hat{G})\) by the formula:
\begin{equation}
  \label{eq:Berezin-Schwartz-kernel}
  \hat{f}(\xi)=\int_G f(x)\,\overline{\chi_\xi(x)}\,dx=
  \int_{\Space{R}{}} f(x)\,e^{-2\pi\rmi \xi x}\,dx.
\end{equation}
This map is unitary on \(\FSpace{L}{2}(G)\).  Pontryagin duality
ensures that the second dual \(\hat{\hat{G}}\) is canonically
isomorphic to \(G\) and provides an expression for the inverse Fourier
transform:
\begin{equation}
  \label{eq:fourier-integral}
  f(x)=
  \int_{\Space{R}{}} \hat{f}(x)\,e^{2\pi\rmi \xi x}\,dx.
\end{equation}
Then, one can interpret the formula~\eqref{eq:PDO-KN-defn} as follows:
\begin{equation}
  \label{eq:AHW-PDO-first}
  [Au](y)=\int_{\widehat{G}} 
  \chi_\xi(y)\,a(y,\xi) \int_{G}u(x)\,\overline{\chi_\xi(x)} \,dx\,d\xi\,,
\end{equation}
where the symbol \(a(x,\xi)\) is a function on \(G\times \hat{G}\).
Since Pontryagin duality and the respective Fourier transforms are
readily available for a locally-compact abelian group, this viewpoint
generates a related theory of PDO on commutative groups,
see~\amscite{RuzhanskyTurunen10a}*{Part~II}.

The situation is different for non-commutative groups. The dual object
\(\hat{G}\) of a non-abelian group \(G\)---the collection
of all equivalence classes of irreducible unitary representations---is
not a group, in general. One can still define the (operator
valued!) Fourier transform by the formula
\begin{displaymath}
  \hat{f}(\xi)=\int_G f(g)\,\xi(g)\, dg, \qquad\text{where } f\in
  \FSpace{L}{1}(G,dg),\quad \xi\in\hat{G}
\end{displaymath}
and \(dg\) is a left-invariant (Haar) measure on \(G\). The inverse
Fourier transform is not as simple as in the commutative case. For
example, on a compact group it is:
\begin{displaymath}
   [\oper{F}^{-1}F] (x) = \sum_{[\xi]\in\hat{G}} \dim(\xi)\,\mathrm{Tr} \left(\xi(x) F (\xi)\right)\,,
\end{displaymath}
where \(\mathrm{Tr}\) denotes the trace of an operator. Thus, for a compact
group an analog of PDO with a symbol  \(a (x, \xi)\) on \(G\times
\hat{G}\) can be defined by, \amscite{RuzhanskyTurunen10a}*{(10.19)}:
\begin{equation}
  \label{eq:oper-dual-gr}
  Af (x) = \sum_{[\xi]\in\hat{G}} \dim(\xi)\,\mathrm{Tr} \left(\xi(x)\, a (x, \xi)\, \hat{f} (\xi)\right).
\end{equation}
Similar formulae were used in the context of the Heisenberg
group~\cite{BahouriFermanian-KammererGallagher12} and other nilpotent
Lie groups~\cite{FischerRuzhansky12a}. Furthermore, the Pontryagin
duality is employed in generalization of Toeplitz
operators%
\index{Toeplitz!operator}%
\index{operator!Toeplitz}~\cite{Mirotin11a} and Wiener--Hopf
factorization%
\index{Wiener--Hopf!factorization}%
\index{factorization, Wiener--Hopf}~\cite{EhrhardtMeeRodmanSpitkovsky07a}, see the cites
papers for details and further references.

\subsection{Covariant Transform}
\label{sec:covariant-transform}

A different approach starts from the observation that operators of
spatial shifts \(f(t) \mapsto f(t-x)\) and operators of
multiplications by exponents \(f(t) \mapsto e^{2\pi\rmi y t} f(t)\)
(i.e. shifts in the frequency space) generate the Schr\"odinger
representation of the non-commutative Heisenberg group
\(\Space{H}{}\)~\cites{Howe80a,Folland89}. As \(C^\infty\)-manifold
\(\Space{H}{}\) can be identified with \(\Space{R}{3}\) and the group
law is:
\begin{equation}
  \label{eq:H-1-group-law}
  \textstyle
  (s,x,y)*(s',x',y')=(s+s'+\frac{1}{2}(xy'-x'y),x+x',y+y').
\end{equation} 
The Schr\"odinger representation \(\Space{H}{}\) is
\begin{equation}
  \label{eq:schroedinger-rep}
  \uir{}{}(s,x,y)f(t)=  e^{\pi \rmi(2s+y(2t-x))}\,  f(t-x).
\end{equation}
We can integrate this representation with the Fourier transform
\(\hat{\sigma}(x,y)\) of a function \(\sigma(q,p)\) on
\(\Space{R}{2}\)~\amscite{Folland89}*{\S~2.1}:
\begin{align}
    [\uir{}{}(\hat{\sigma})f](t)=&\int_{\Space{R}{2}}
    \hat{\sigma}(x,y)\,\uir{}{}(0,x,y)f(t)\,dx\,dy \nonumber \\
    =& \int_{\Space{R}{2}} \hat{\sigma}(x,y)\,e^{\pi \rmi y(2t-x)}\,
    f(t-x)\,dx\,dy \nonumber \\
    =& \int_{\Space{R}{2}} \textstyle \sigma(q,\frac{1}{2}(r+t))\,
    e^{2\pi\rmi q(r-t)} \, f(r)\,dq\,dr.
  \label{eq:pdo-def}
\end{align}
Up to different letters, this is exactly PDO~\eqref{eq:PDO-Weyl-defn}
with the Weyl symbol \(\sigma\).

It may not be obvious that an introduction of the non-commutative
Heisenberg group produces any advantage over commutative
Pontryagin duality. Probably, it explains why this direction, rooted in
Weyl's original works and spectacularly developed
in~\cites{Howe80a,Howe80b}, was not widely adopted (see, however, 
remarkable
exceptions~\cites{Dynin75,Dynin76,Ratcliff85a,Folland89,Folland94}).

Benefits, which can be challenged within the classical PDOs, become
more explicit when we move to a general setup. A transition to a
non-commutative underlining group does not become an issue since
non-commutativity is already in the scheme. Thus, the construction of
non-abelian PDOs is different in a computational sense rather than
conceptually. Moreover, historically PDOs appeared as ``SIOs varying
from point to point'' and these roots were preserved
in~\cite{Dynin75,Dynin76}.  We will recover them in examples with the
Dynin group below.

Such a development is very straightforward for nilpotent Lie groups,
as was already hinted in~\cite{Dynin75,Dynin76,Howe80b}. Thus, the
concept of relative convolutions%
\index{relative convolution}%
\index{convolution!relative}~\cite{Kisil94e} was initially
developed in the nilpotent setting.  However, the approach is also
usable for non-compact non-commutative non-exponential (e.g.
semisimple) Lie groups as well.  The present paper provides a brief
illustration to this claim.

\section{Groups and Representations}
\label{sec:preliminaries}
Our construction is based on groups and representation theory. It is
connected to the covariant
transform~\cites{Kisil11c,Kisil10c,Kisil09d,Kisil12d}, which
consolidates a large collection of results linked to
wavelets/coherent states~\cites{Kisil94e,Kisil98a}.

\subsection{Main Examples}
\label{sec:main-examples}

For the sake of brevity, we explicate our approach only by the
following four examples. However, they do not exhaust all possible
applications. 

\subsubsection{The Heisenberg Group}
\label{sec:heisenberg-group-1}
For simplicity, we use only the smallest one-dim\-en\-sio\-nal
\emph{Heisenberg group}%
\index{Heisenberg!group}%
\index{group!Heisenberg}
\(\Space{H}{}\) consisting of points \((s,x,y)\in
\Space{R}{3}\)~\cites{Folland89,Howe80b}. The group law on
\(\Space{H}{}\) is given by~\eqref{eq:H-1-group-law}.  The Heisenberg
group is a non-commutative nilpotent Lie group with the centre
\begin{displaymath}
  Z=\{(s,0,0)\in \Space{H}{}, \ s \in \Space{R}{}\}.
\end{displaymath}
The Lie algebra \(\algebra{h}\) is realised by the following
left-(right-)invariant vector fields:
\begin{equation}
  \label{eq:h-lie-algebra}
  \textstyle  S^{l(r)}=\pm{\partial_s}, \quad
  X^{l(r)}=\pm\partial_{x}-\frac{1}{2}y{\partial_s},  \quad
  Y^{l(r)}=\pm\partial_{y}+\frac{1}{2}x{\partial_s}.
\end{equation}
They satisfy to the \emph{Heisenberg commutator relations}%
\index{Heisenberg!commutator relations}%
\index{relations!Heisenberg commutator} \([X,Y]=S\) and \([X,S]=[Y,S]=0\).





\subsubsection{Abstract Heisenberg--Weyl (AHW) Group}
\label{sec:abstr-heis-weyl}

Let \(G\) be a locally compact abelian group. Pontryagin duality tells
that the collection \(\hat{G}\) of all unitary characters of \(G\) is
a locally compact group as well. For
example~\amscite{KirGvi82}*{\S~IV.2.1}, \(\hat{\Space{R}{}}=\Space{R}{}\),
\(\hat{\Space{Z}{}}=\Space{T}{}\), \(\hat{\Space{T}{}}=\Space{Z}{}\),
where \(\Space{T}{}\) is the group of unimodular complex numbers. The
group operations on both \(G\) and \(\hat{G}\) are denoted by \(+\) and
their units are written as \(0\).

We form a new group \(\tilde{G}\) as the set
\(\Space{T}{}\times G\times \hat{G}\) with the group
law~\citelist{\cite{Mackey49a} \cite{Mackey70a} \cite{Folland89}*{\S~1.11}}:
\begin{displaymath}
  (z_1,g_1,\chi_1)*(z_2,g_2,\chi_2)=(z_1 z_2 \chi_2(g_1)
  ,g_1+g_2,\chi_1+\chi_2),
\end{displaymath}
where \(z_i\in \Space{T}{}\), \(g_i\in G\), \(\chi_i\in \hat{G}\),
\(i=1,2\).  In general, \(\tilde{G}\) is a non-commutative locally
compact group. The unit is \((1,0,0)\) and the inverse of
\((z,g,\chi)\) is \((\bar{z}{\chi}(g),-g,-\chi)\). The centre of
\(\tilde{G}\) consists of elements \((z,0,0)\), \(z\in\Space{T}{}\).
The left- and right-invariant measures coincide with the product of
invariant measures of \(\Space{T}{}\), \(G\) and \(\hat{G}\). 

For \(G=\Space{R}{}\), the group \(\tilde{G}\) is the
polarised Heisenberg group~\cite{Folland89}*{\S~1.2} in the reduced
form~\cite{Folland89}*{\S~1.3}. Thus, for a general \(G\), we call
\(\tilde{G}\) the \emph{abstract Heisenberg--Weyl (AHW) group}%
\index{abstract Heisenberg--Weyl (AHW) group}%
\index{group!abstract Heisenberg--Weyl (AHW)}%
\index{Heisenberg--Weyl! abstract group (AHW)}%
\index{AHW group|see{abstract Heisenberg--Weyl group}}. Another basic example
of the AHW group is
\(\tilde{\Space{T}{}}=\Space{T}{}\times\Space{T}{}\times\Space{Z}{}\),
with the group law:
\begin{displaymath}
  (z_1,w_1,k_1)*(z_2,w_2,k_2)=(z_1z_2 w_1^{k_2}, w_1 w_2, k_1+k_2),
\end{displaymath}
where \(z_i, w_i\in\Space{T}{}\) and \(k_i\in\Space{Z}{}\), \(i=1,2\)
and the operation on \(\Space{T}{}\) written as multiplication of complex numbers.
Of course, by the Pontryagin duality \(\widetilde{\hat{G}}\) is
isomorphic to \(\tilde{G}\), in particular,
\(\tilde{\Space{T}{}}\) is isomorphic to \(\tilde{\Space{Z}{}}\). Our
consideration of \(\tilde{G}\) shall be compared
with~\cite{Wildberger05b}.

\subsubsection{The Dynin Group}
\label{sec:dynin-group}

Extending the Heisenberg group, consider a Lie algebra \(\algebra{d}\)
spanned by the basis \(\{Z,T,U,V,S,X,Y\}\) defined by the following
non-vanishing commutators~\cite{Dynin75}:
\begin{align}
  [X,Y]&=S,&  [X,U]&=Z,&  [Y,V]&=Z, \\
   [S, T]&=Z,& [X,T]&\textstyle=-\frac{1}{2} V,&  [Y,T]&\textstyle=\frac{1}{2}U.
\end{align}
Thus, the Lie algebra \(\algebra{d}\) is nilpotent step \(3\). It is
generated by its elements \(X\), \(Y\) and \(T\) and their
commutators. The multiplication on a group \(\Space{D}{}\), obtained from
\(\algebra{d}\) by exponentiation, is:
\begin{align}
  \label{eq:dynin-group-law}
  \lefteqn{(z,t,u,v,s,x,y)*(z',t',u',v',s',x',y')}\quad\qquad
  \\
    =&\textstyle
    (z+z'+\frac{1}{2}(st'-s't)+\frac{1}{2}(xu'-x'u)+\frac{1}{2}(yv'-y'v)\nonumber \\
    &\textstyle\qquad+\frac{1}{24}(yx't-xy't+y'xt'-x'yt'),\nonumber \\
    &\ \textstyle t+t',\, u+u'+\frac{1}{4}(yt'-y't),\,
    v+v'-\frac{1}{4}(xt'-x't),\nonumber \\
    &\ \textstyle s+s'+\frac{1}{2}(xy'-x'y),x+x',y+y').
    \nonumber
\end{align} 

A representation \(dR\) of \(\algebra{d}\) appears if we extend the
representation of the Lie algebra \(\algebra{h}\) by the
left-invariant vector fields~\eqref{eq:h-lie-algebra} with the
representation of the additional operators \(U\), \(V\), \(T\), \(Z\)
as operators of multiplication:
\begin{equation}
  \label{eq:dynin-mult-rep}
  dR^{U}=xI,\qquad dR^V=yI, \qquad dR^T=sI,\qquad dR^Z=I. 
\end{equation}
This representation is connected with the algebra generated by
convolutions on the Heisenberg group and operators of multiplications
by functions, see Example~\ref{it:dynin-source} below. The group \(\Space{D}{}\) was
used in papers~\cites{Dynin75,Dynin76}, thus we call it the
\emph{Dynin group}.  It is a special (but, probably,
the most important) example of meta-Heisenberg group~\cite{Folland94}.
It is also a subgroup of the group studied in~\cite{Ratcliff85a}.



\subsubsection{The group $\SU$}
\label{sec:group-su}

The group $\SU$~\citelist{\cite{Lang85}*{\S~IX.1}
  \cite{MTaylor86}*{\S~8.1}} consists of \(2\times 2\) matrices with
complex entries of the form \(
\begin{pmatrix}
  \alpha &\beta \\\bar{\beta}&\bar{\alpha}
\end{pmatrix}\) and unit determinant:
\(\modulus{\alpha}^2-\modulus{\beta}^2=1\). The multiplication is
given by matrix multiplication and is not commutative. The maximal compact
subgroup \(K\) of diagonal matrices \(
\begin{pmatrix}
  e^{\rmi\phi}&0\\0&e^{-\rmi\phi}
\end{pmatrix}\) is isomorphic to the unit circle. Its presence
indicates that the group is not exponential, that is the exponent map
from the Lie algebra to the group is not a bijection. This is group is
not compact.

The group \(\SU\) acts by M\"obius transformations of the unit disk
(Example~\ref{it:SU-hs-disk} below) and is very important in complex
analysis, cf.~\cite{Coburn12a}.  Intimate connections of other
subgroups of \(\SU\) with the hypercomplex numbers is described
in~\cites{Kisil12a,Kisil11c,Kisil09e}, we do not touch this interesting
topic in this paper.

\subsection{Induced Representations}
\label{sec:induc-repr}

The general scheme of induced representations is as follows, see
also~\citelist{\amscite{Kirillov76}*{\S~13.2}
  \amscite{MTaylor86}*{Ch.~5} \amscite{Folland95}*{Ch.~6}
  \amscite{Kisil97c}*{\S~3.1} \amscite{Kirillov04a}*{\S~V.2} \cite{Mackey70a}}.  

Let \(G\) be a locally compact group and let \(H\) be its subgroup.  Let \(X=H
\backslash G\) be the corresponding right coset space and \(\map{s}: X
\rightarrow G\)%
\index{map!$\map{s}$}\index{$s$@$\map{s}$
  (section)}\index{section!$\map{s}$} %
be a continuous function (section)~\amscite{Kirillov76}*{\S~13.2} which
is a right inverse to the natural projection \(\map{p}:G\rightarrow H
\backslash G\).%
\index{$p$@$\map{p}$ map}%
\index{map!$\map{p}$} Then, any \(g\in G\) has a unique decomposition
of the form \(g=h * \map{s}(x)\) where \(x=\map{p}(g)\in X\) and
\(h\in H\). We define the map \(\map{r}: G\rightarrow H\):%
\index{map!$r$}\index{$r$ (map)}
\begin{equation}
  \label{eq:r-map1}
  \map{r}(g)=g *{\map{s}(x)}^{-1}, \qquad \text{ where } x=\map{p}(g).
\end{equation}
Note, that \(X\) is a right homogeneous space with the \(G\)-action
defined in terms of \(\map{p}\) and \(\map{s}\) as follows:
\begin{equation}
  \label{eq:g-action}
  g: x  \mapsto x\cdot g=\map{p}(\map{s}(x)*g),
\end{equation}
where \(*\) is the multiplication on \(G\).

\begin{example}\upshape
  \begin{enumerate}
  \item \label{it:heisenberg-map}
   For the Heisenberg group \(\Space{H}{}\) we can consider the
    subgroup \(Z=\{(s,0,0)\such s\in\Space{R}{}\}\). The corresponding
    homogeneous space is \(Z\backslash \Space{H}{}=\{(0,x,y)\such
    (x,y)\in\Space{R}{2n}\}\). Using the maps \(\map{p}:
    (s',x',y')\mapsto (x',y')\) and \(\map{s}:(x',y')\mapsto(0,x',y')\)
    we calculate the action:
    \begin{equation}
      \label{eq:H1-shift-plane}
      (s,x,y): (x',y') \mapsto (x+x',y+y').
    \end{equation}
    There is also a subgroup
    \begin{equation}
      \label{eq:H-x-subgroup}
      H_x=\{(s,0,y)\in\Space{H}{} \such     s,y\in\Space{R}{}\}
    \end{equation}
    and the respective homogeneous space is parametrised by the real
    line. Using the maps \(\map{p}: (s',x',y')\mapsto x'\) and
    \(\map{s}: x' \mapsto (0,x',0)\) we find the action of
    \(\Space{H}{}\) on \(H_x\backslash \Space{H}{}\):
    \begin{equation}
      \label{eq:H1-act-real-line}
      (s,x,y): x' \mapsto x+x'.
    \end{equation}
  \item For an AHW group%
\index{abstract Heisenberg--Weyl (AHW) group}%
\index{group!abstract Heisenberg--Weyl (AHW)}%
\index{Heisenberg--Weyl! abstract group (AHW)} \(\tilde{G}\), there are also two commutative
    subgroups: the centre \(Z=\{(z,0,0)\such z\in\Space{T}{}\}\) and 
    \begin{equation}
      \label{eq:AHW-G-subgroup}
      H_G=\{(z,0,\chi)\in\tilde{G} \such     z\in\Space{T}{}, \ \chi\in\hat{G}\}.
    \end{equation}
    The natural maps \(\map{s}\) and respective actions on the homogeneous
    spaces are similar to the above particular case of \(\Space{H}{}\sim\tilde{\Space{R}{}}\).
  \item For the Dynin group \(\Space{D}{}\), consider the commutative
    subgroup
    \begin{equation}
      \label{eq:M-subgroup-Dynin}
      M=\{(z,t,u,v,0,0,0)\in \Space{D}{} \such (z,t,u,v)\in\Space{R}{4}\}.
    \end{equation}
    The homogeneous space \(M\backslash\Space{D}{}\) can be identified
    with \(\Space{H}{}\) through the map
    \(\map{s}(s,x,y)=(0,0,0,0,s,x,y)\). The corresponding action is
    in the essence the group law~\eqref{eq:H-1-group-law} of
    \(\Space{H}{}\):
    \begin{equation}
      \label{eq:Dynin-Heisenberg-shift}
      (z,t,u,v,s,x,y)\cdot(s',x',y')=(s+s'+\frac{1}{2}(xy'-x'y),x+x',y+y').
    \end{equation}
    
  \item 
    \label{it:SU-hs-disk}
    For the group \(\SU\) and its subgroup \(K\) we identify
    \(K\backslash \SU\) with the open unit disk%
    \index{unit disk}%
    \index{disk, unit} \(D=\{z\in\Space{C}{} \such \modulus{z}<1 \}\). Defining maps:
    \begin{equation}
      \label{eq:SU-maps}
      \map{p}:
      \begin{pmatrix}
        \alpha&\beta\\\bar{\beta}&\bar{\alpha}
      \end{pmatrix}\mapsto \frac{\beta}{{\alpha}}
      \quad\text{ and }\quad
      \map{s}: z\mapsto
      \frac{1}{\sqrt{1-\modulus{z}^2}}\begin{pmatrix}
        1&z\\\bar{z}&1
      \end{pmatrix}, \quad \text{ where } \modulus{z}<1
    \end{equation}
    we deduce the respective action:
    \begin{displaymath}
      z\cdot g=\frac{\bar{\alpha}z+\beta }{\bar{\beta} z+\alpha  },
      \qquad \text{ where } g=\begin{pmatrix}
        \alpha&\beta\\
        \bar{\beta}&\bar{\alpha}
      \end{pmatrix}.
    \end{displaymath}
    This is a linear-fractional (M\"obius) transformation of the unit
    disk~\cite{Kisil12a}*{Ch.~10}.
  \end{enumerate}
\end{example}

For \(G\) and \(H\), we respectively denote the right Haar measures
\(dg\) and \(dh\), the corresponding modular functions are \(\Delta_G\) and
\(\Delta_H\).  Then, there is a measure \(dx\) on \(X=H\backslash G\)
defined up to a scalar factor by the
identities~\amscite{Kirillov76}*{\S~9.1(5\('\)-6\('\))}:
\begin{equation}
  \label{eq:hs-measure}
  dg=\frac{\Delta_G(h)}{\Delta_H(h)}\,dx\,dh, \qquad \text{ where }
  g=h\map{s}(x) . 
\end{equation}
The measure \(dx\) transforms under \(G\) action~\eqref{eq:g-action}
by, see~\amscite{Kirillov76}*{\S~9.1(7\('\)-8\('\))}: 
\begin{equation}
  \label{eq:hs-measure-transf}
  \frac{d(x\cdot    g)}{dx}
  =\frac{\Delta_H(h(x,g))}{\Delta_G(h(x,g))},\qquad
  \text{ where} \quad \map{s}(x)g=h(x,g)\map{s}(x\cdot g).
\end{equation}

In many cases, e.g. for all nilpotent Lie
groups~\amscite{Kirillov04a}*{\S~3.3.2}, unitary representations are
induced by characters---one dimensional linear representations---of
its subgroups. Thus we present here the induction from
a character only. Let \(\chi: H \rightarrow \Space{C}{}\) be a unitary
character of \(H\). Consider the space of functions on \(G\) with the
property:
\begin{equation}
  \label{eq:induced-space-prop}
  F(hg)=\chi(h) F(g),
\end{equation}
the space is obviously invariant under right translations.  The restriction of the
right regular representation:%
\index{regular!representation!right}%
\index{representation!regular!right}%
\index{right!regular representation}
\(R(g):  f(g')\mapsto f(g'g)\) to this space is called an \emph{induced
  representation}%
\index{induced!representation}%
\index{representation!induced} in the sense of
Mackey~\citelist{\amscite{Kirillov76}*{\S~13.2}
  \amscite{Kirillov04a}*{\S~V.2}}.
  
Consider the \emph{lifting}\index{lifting} \(\oper{L}_\chi:
\FSpace{C}{b}(X)\rightarrow \FSpace{C}{b}(G)\) of continuous bounded
functions:
\begin{equation}
  \label{eq:lifting-defn}
  F(g)=[\oper{L}_\chi f](g)=\chi(h) f(\map{p}(g)),\qquad
  f(x)\in \FSpace{C}{b}(X).
\end{equation}
The function \(F(g)\) has the property~\eqref{eq:induced-space-prop}.
The same expression~\eqref{eq:lifting-defn} defines a bijection from
\(\FSpace{L}{p}(X)\) to certain space \(\FSpace[\chi]{L}{p}(G)\),
which is invariant under right translations.  A right inverse
map---the \emph{pulling}\index{pulling}---\(\oper{P}: \FSpace[\chi]{L}{p}(G)
\rightarrow \FSpace{L}{p}(X)\) is defined by:
\begin{equation}
  \label{eq:pulling-defn}
  f(x)=[\oper{P}F](x)=F(\map{s}(x)),\qquad F(g)\in \FSpace[\chi]{L}{p}(G).
\end{equation}
The norm on \(\FSpace[\chi]{L}{p}(G)\) is introduced in such a way
that both the lifting\index{lifting} and pulling are isometries.

Since \(\FSpace[\chi]{L}{p}(G)\) is invariant under the right shifts,
lifting and pulling intetwine the restriction
\(R|_{\FSpace[\chi]{L}{p}(G)}\) of the right regular representation%
  \index{regular!representations}%
  \index{representations!regular}
\(R\) with the representation \(\uir{}{\chi}(g)= \oper{P} \circ R(g)
\circ \oper{L}_\chi\).  It is the second form of the induced
  representation%
\index{induced!representation}%
\index{representation!induced}. Its
realisation \(\uir{}{\chi}\) in a space of complex-valued functions on
\(X\), cf.~\amscite{Kirillov76}*{\S~13.2(7)--(9)} is:
\begin{equation}
  \label{eq:def-ind}
  [\uir{}{\chi}(g) f](x)= \chi(\map{r}(\map{s}(x)*g))\,
  f(x\cdot g),
\end{equation}
where \(g\in G\), \(x\in X\), \(h\in H\) and \(\map{r}: G \rightarrow
H\), \(\map{s}: X \rightarrow G\) are maps defined above;
\(*\)~denotes multiplication on \(G\) and \(x\cdot g\) denotes the
action~\eqref{eq:g-action} of \(G\) on \(X\) from the right.

For the case of an unimodular group \(G\) and an unimodular subgroup
\(H\subset G\) (which is automatic for a nilpotent \(G\)), the
representation~\eqref{eq:def-ind} is unitary in \(\FSpace{L}{2}(X)\).
In the case of a non-unimodular (sub)group, we need an additional
factor
\(\left[\frac{\Delta_H(h(x,g))}{\Delta_G(h(x,g))}\right]^{\frac{1}{2}}\)
to make \(\uir{}{\chi}\) unitary, cf.~\eqref{eq:hs-measure-transf}
and~\amscite{Kirillov76}*{\S~13.2(3)}. 

\begin{example}\upshape
  \begin{enumerate}
  \item For the centre \(Z\) of \(\Space{H}{}\), the map \(\map{r}:
    \Space{H}{}\rightarrow Z\) is \(\map{r}(s,x,y)=(s,0,0)\). The
    character \(\chi_\myhbar(s,0,0)=e^{2\pi\rmi\myhbar s}\) of \(Z\)
    together with the action~\eqref{eq:H1-shift-plane} produces the
    unitary Fock--Segal--Bargmann (FSB) representation%
    \index{FSB!representation}%
    \index{representation!FSB}%
    \index{Fock--Segal--Bargmann!representation|see{FSB representation}}%
    \index{representation!Fock--Segal--Bargmann|see{FSB representation}}~\citelist{\amscite{Folland89}*{\S~1.6}
      \cite{Kisil10a}}:
    \begin{equation}
      \label{eq:H1-rep-plane}
      [\uir{F}{\myhbar}(s,x,y)f] (x',y') = e^{\pi \rmi \myhbar  (2s+x'y-xy')} f(x'+x,y'+y).
    \end{equation}
    We identify a point \((x,y)\) of the homogeneous space
    \(Z\backslash \Space{H}{}\) with the complex number \(z=x+\rmi
    y\). Then, the representation~\eqref{eq:H1-rep-plane} can be
    stated in the complex form:
    \begin{equation}
      \label{eq:H1-rep-plane-compl}
      [\uir{F}{\myhbar}(\map{s}(z))f] (z') = e^{\pi\myhbar  (z\bar{z}'-\bar{z}z')/2} f(z+z').
    \end{equation}

    For the subgroup \(H_x\)~\eqref{eq:H-x-subgroup}, the map
    \(\map{r}(s,x,y)=(s+\frac{1}{2}xy,0,y)\). A character
    \(\chi_\myhbar(s,0,y)=e^{2\pi\rmi\myhbar s}\) and the
    action~\eqref{eq:H1-act-real-line} produce the Shr\"odinger
    representation%
    \index{Shr\"odinger!representation}%
    \index{representation!Shr\"odinger}~\citelist{\cite{Folland89}*{\S~1.3}  \cite{Kisil10a}}:
    \begin{equation}
      \label{eq:Schroedinger-alt}
      [\uir{}{\myhbar}(s,x,y)f] (x')= e^{\pi\rmi\myhbar (2s+2x'y+xy)} f(x'+x).
    \end{equation}
    Clearly, \(\uir{}{1}(-s,-x,-y)\) coincides with the
    representation~\eqref{eq:schroedinger-rep}. Furthermore, it is
    known that they are unitary equivalent to FSB
    representation~\eqref{eq:H1-rep-plane}.
  \item For an AHW group%
\index{abstract Heisenberg--Weyl (AHW) group}%
\index{group!abstract Heisenberg--Weyl (AHW)}%
\index{Heisenberg--Weyl! abstract group (AHW)} \(\tilde{G}\), we proceed in a similar
    fashion. The  character \(\nu(z,0,0)=z^k\) of the
    centre induce the representation on \(\FSpace{L}{}(G\times\hat{G})\):
    \begin{equation}
      \label{eq:AHW-rep-FSB}
      [\uir{F}{k}(z,g,\chi) f](g',\chi')=(z\chi(g'))^k\, f(g+g',\chi+\chi').
    \end{equation}
    For the subgroup \(H_G\)~\eqref{eq:AHW-G-subgroup}, the map
    \(\map{r}(z,g,\chi)=(z,0,\chi)\) and the character
    \(\nu(z,0,\chi)=z^k\) induces the representation on
    \(\FSpace{L}{}(G)\):
    \begin{equation}
      \label{eq:AHW-Shroedinger}
      [\uir{}{k}(z,g,\chi) f](g')=(z\chi(g'))^k\, f(g+g').
    \end{equation}
    The classification of irreducible representations of AHW group was
    provided in~\cite{Mackey49a} in a way which generalised the
    Stone--von Neumann theorem for \(\Space{H}{n}\), see
    also~\cite{Mackey70a}. 
  \item \label{it:Dynin-represent}
    For the subgroup \(M\)~\eqref{eq:M-subgroup-Dynin} of
    \(\Space{D}{}\), the map
    \begin{displaymath}
      \map{r}
      (z,t,u,v,s,x,y) =\textstyle (z+\frac{1}{2}st+\frac{1}{2}xu+\frac{1}{2}yv,\,
      \textstyle t,\, u+\frac{1}{4}yt,\, v-\frac{1}{4}xt,\,
      \textstyle 0,\,0,\,0).
    \end{displaymath}
    The character \(\chi(z,t,u,v,0,0,0)=e^{2\pi\rmi\myh z}\) of \(M\)
    induces the representation:
    \begin{align}
      \label{eq:Dynin-represent}
        \lefteqn{[\uir{}{\myh}(z,t,u,v,s,x,y)f](s',x',y')}\qquad \\
        =& e^{\pi\rmi\myh(2z+2s't+st+\frac{1}{4}(x'y-xy')t
          +(2x'+x)u +(2y'+y)v)}  \nonumber\\
        & {} \times \textstyle
        f(s+s'+\frac{1}{2}(x'y-xy'),x+x',y+y').\nonumber
    \end{align}
      In particular, the action of the subgroup \(M\) reduces to multiplication:
      \begin{equation}
        \label{eq:Dynin-represent-mult}
        [\uir{}{\myh}(z,t,u,v,0,0,0)f](s',x',y')=e^{2\myh\pi\rmi(z +ts'+ux'+vy')}f(s',x',y').
    \end{equation}
    On the other hand, the operator \(\uir{}{\myh}(0,0,0,0,s,x,y)\) is
    the shift~\eqref{eq:Dynin-Heisenberg-shift} by \((s,x,y)\) on
    \(\Space {H}{}\). The corresponding infinitesimal actions
    are~\eqref{eq:h-lie-algebra} and~\eqref{eq:dynin-mult-rep}. This
    representation was used in~\cite{Dynin75}, see
    also~\cite{Folland94}*{\S~3}.

  \item \label{it:SU-representation}
    For \(G=\SU\) and \(H=K\) we calculate \(\map{r}
    \begin{pmatrix}
      \alpha&\beta \\
      \bar{\beta}&\bar{\alpha}
    \end{pmatrix}=
    \begin{pmatrix}
      \frac{\alpha}{\modulus{\alpha}}&0\\
      0&\frac{\bar{\alpha}}{\modulus{\alpha}}
    \end{pmatrix}\). Let a character \(\chi\) of \(K\) be \(\chi
    \begin{pmatrix}
      e^{\rmi \phi}&0\\0&e^{-\rmi \phi}
    \end{pmatrix}=e^{-2\rmi \phi}\), then the induced representation
    acts on \(\FSpace{L}{2}(D)\) as follows:
    \begin{equation}
      \label{eq:SU-action}
      [\uir{}{}(g) f](z) = \frac{\modulus{\bar{\beta}z+\alpha }^2}{(\bar{\beta}z+\alpha)^2 }\, 
      f\left(\frac{\bar{\alpha}z+\beta }{\bar{\beta} z+\alpha
        }\right)
      = \frac{\beta\bar{z}+\bar{\alpha }}{\bar{\beta}z+\alpha}\, 
      f\left(\frac{\bar{\alpha}z+\beta }{\bar{\beta} z+\alpha  }\right).
    \end{equation}
    Since we are not in the unimodular setting now, we calculate the
    invariant measure on the unit disk to be
    \(\left.(1-\modulus{z}^2)\right.^{-2}\,dz\wedge d\bar{z}\). The
    representation~\eqref{eq:SU-action} is unitary and belongs to the
    discrete series~\amscite{Lang85}*{\S~IX.3}. In contrast to
    equivalent representations used in complex analysis, our
    expression~\eqref{eq:SU-action} has clearer composition formula,
    cf.~\cite{Coburn12a}*{(**)}. However, \eqref{eq:SU-action} does
    not preserves usual analyticity. We can use either
    conformal-invariant modification of the Cauchy--Riemann
    equations~\cite{Kisil11c}*{\S~5.3}, or introduce an additional
    peeling map, which intertwines our representation with the more
    common one, acting in the space of analytic functions.
  \end{enumerate}
\end{example}

\section{Covariant Transform}
\label{sec:covariant-transform-1}
Representation theory is behind many important calculations in
analysis, this is illustrated in the present section. The
group-theoretical foundations of coherent states/wavelets are
well-known and widely
appreciated~\cites{AliAntGaz00,Perelomov86,Kisil02a,%
Kisil95i,Kisil12b,Kisil11c,Kisil98a}.

\subsection{Induced Covariant Transform}
\label{sec:induc-covar-transf}
The following definition is a general template, which admits various
specialisations adjusted to particular cases. 
\begin{defn} \cite{Kisil09d} Let \(\uir{}{}\) be a representation of
  \(G\) in a vector space \(V\). For a vector space \(U\) and an
  operator \(F: V \rightarrow U\), the \emph{covariant transform}%
  \index{covariant!transform}%
  \index{covariant!transform}%
  \index{transform!covariant} is the map:
  \begin{equation}
    \label{eq:covariant-trans}
    [\oper{W}_F v](g)=F(\uir{}{}(g)v),\qquad v\in V,
    \quad g\in   G,
  \end{equation}
  to \(U\)-valued functions on \(G\). In this context we call \(F\) a
  \emph{fiducial operator}%
  \index{fiducial operator}%
  \index{operator!fiducial}.
\end{defn}

An important particular case of the above definition is provided by a
linear functional \(F\in V^*\), the covariant transform produces
\emph{matrix coefficients}%
\index{matrix coefficient!representation, of}%
\index{representation!matrix coefficient} of the
representation~\cite{HoweTan92}*{Ex.~I.1.2.12}. In the case of a Hilbert
space \(V\), such a functional is provided by a pairing with a vector
\(f\in V\), which is known as \emph{mother wavelet}%
\index{wavelet!mother}%
\index{mother wavelet} or \emph{vacuum state}%
\index{vacuum state}%
\index{state!vacuum}~\cites{AliAntGaz00,Perelomov86}.  Then, the
covariant transform becomes the wavelet transform:
\begin{equation}
  \label{eq:wavelet-transform}
   \tilde{v}(g)=[\oper{W}_f
   v](g):=\scalar{\uir{}{}(g)v}{f}=\scalar{v}{\uir{*}{}(g)f}=\scalar{v}{f_g}, 
\end{equation}
where \(f_g=\uir{*}{}(g)f\) are called wavelets.  The image
\(\oper{W}_f v\) is a scalar-valued function. The scalar case is very
important, however, it does not cover all interesting situations,
see~\cite{Kisil12b} and Example~\ref{ex:wick-symbol}. We also may
require a functional \(F\) associated to a singular mother wavelet,%
\index{wavelet!mother!singular}%
\index{mother wavelet!singular}%
\index{singular!mother wavelet}
i.e. a distribution, cf.~\cite{Kisil98a}*{\S~2.3}.

If the representation \(\uir{}{}\) and the operator \(F\) are bounded,
then the image of \(\oper{W}_F\) consists of bounded functions on
\(G\). Weak continuity of \(\uir{}{}\) suffices for continuity of
\(\oper{W}_F v\). An important property of \(\oper{W}_F\) is as
follows.
\begin{lem}
  \label{le:covariant-intertwining}
  The covariant transforms intertwines the left%
  \index{left regular representations!}%
  \index{regular!representations!left}%
  \index{representations!regular!left} \(\Lambda(g):
  f(g')\mapsto f(g^{-1}g')\) and right \(R(g)\) regular
  representations%
  \index{regular!representation!right}%
  \index{representation!regular!right}%
  \index{right!regular representation} of \(G\) with the following
  actions of \(\uir{}{}\):
  \begin{equation}
    \label{eq:covariant-intertwine}
    R (g) \oper{W}_F = \oper{W}_F \uir{}{}(g)
    \quad\text{ and }\quad
     \Lambda (g) \oper{W}_F= \oper{W}_{ F\circ \uir{}{}(g^{-1})}
    \quad\text{ for all } g\in G.
  \end{equation}
\end{lem}
There is the following simple but useful consequence of the above Lemma.
\begin{cor} \cite{Kisil11c}*{Cor.~5.8}
  \label{co:cauchy-riemann-integ}
  Let \(\uir{}{}\) be a linear representation
   of a group \(G\) on a space \(V\), which has an adjoint representation
  \(\uir{*}{}\) on the dual space \(V^*\). Let a mother wavelet \(f\in V^*\)  satisfy the equation 
  \begin{displaymath}
    \int_{G} a(g)\, \uir{*}{}(g) f\,dg=0,
  \end{displaymath}
  for a fixed distribution \(a(g) \) and a (not necessarily
  invariant) measure \(dg\). Then,  any wavelet transform
  \(\tilde{v}=\scalar{v}{\uir{*}{}(g)f}\) obeys the following
  right-invariant condition:
  \begin{equation}
    \label{eq:cauchy-riemann-op-abstract}
    D \tilde{v}=0,\qquad \text{where} \quad D=\int_{G} \bar{a}(g)\, \Lambda (g) \,dg,
  \end{equation}
  with \(\Lambda \) being the left regular representation of \(G\).
\end{cor}
As we will see below, the above distribution \(a\) is often a
linear combination of derivatives of the Dirac's delta functions, therefore
the operator \(D\) turns to be a differential operator. Further
examples can be found in~\cite{Kisil12d}*{Ex.~5.9--11}.

Often we need only a part of covariant transform. For a Lie group 
\(G\) and its subgroup \(H\), we fix a continuous section \(\map{s}:
H\backslash G \rightarrow G\), which is a right inverse to the
projection \(\map{p}: G \rightarrow H\backslash G\).
\begin{defn} 
  \cite{Kisil11c}*{\S~5.1} Let \(F:V\rightarrow U\) intertwine the restriction of
  \(\uir{}{}\) to \(H\)  with a character \(\chi\) of \(H\):
  \(F(\uir{}{}(h) v)=\chi(h)F(v)\) for 
  all \(h\in H\), \(v\in V\). Then, the
  \emph{induced covariant transform}%
  \index{induced!covariant transform}%
  \index{covariant!transform!induced}%
  \index{transform!covariant!induced} is:
  \begin{equation}
    \label{eq:induced-covariant-trans}
    [\oper{W}_F v](x)=F(\uir{}{}(\map{s}(x))v),\qquad v\in V,
    \quad x\in    H\backslash G.
  \end{equation}
\end{defn}
Under our assumptions, the induced covariant transform intertwines
\(\uir{}{}\) with the representation induced from \(H\) by the
character \(\chi\). To use the
condition~\eqref{eq:cauchy-riemann-op-abstract} for the induced
covariant transform, we need first apply the lifting\index{lifting}
\(\oper{L}_\chi\)~\eqref{eq:lifting-defn} to \(\oper{W}_F v\) and then
the operator \(D\). A collection of such
conditions~\eqref{eq:cauchy-riemann-op-abstract} can characterise the
image \(\oper{W}_f V\) among all functions on \(X\), see
Examples~\ref{it:FSB-reproducing-space} and~\ref{it:SU-Bergman-proj}
below.  

In many cases, e.g. for square integrable representations%
\index{square integrable!representations}%
\index{representations!square integrable} and an admissible mother
wavelet \(v\in V\), the image space of the covariant transform is a
reproducing kernel Hilbert space~\cite{AliAntGaz00}*{Thm.~8.1.3}. That
means that for any function \(v \in \oper{W}_f V\) we have the
integral reproducing formula:
\begin{equation}
  \label{eq:reproducing-integral}
  v(y)=\int_X v(x)\,\bar{k}_y(x)\,dx,
\end{equation}
where the \emph{reproducing kernel}%
\index{reproducing kernel}%
\index{kernel!reproducing} \(k_y\) provides the twisted convolution
with the normalised covariant transform \(\oper{W}_f
(\uir{}{}(\map{s}(y)^{-1}) f)\) for the mother wavelet \(f\), see
Cor.~\ref{co:reproducing-formula}. For a function \(v\not\in
\oper{W}_f V\), the right-hand side of~\eqref{eq:reproducing-integral}
defines its projection to the space \(\oper{W}_f V\).

\begin{example}
  \begin{enumerate}
  \item \label{it:FSB-reproducing-space}
    For \(G=\Space{H}{}\), \(H=Z\)  and the representation
    \(\uir{}{\myhbar}\)~\eqref{eq:Schroedinger-alt} on
    \(\FSpace{L}{2}(\Space{R}{})\), we have
    \(\uir{}{ \myhbar}(s,0,0)=e^{2\pi\rmi \myhbar s}\). Thus any function \(f
    \in \FSpace{L}{2}(\Space{R}{})\) is suitable for the induced
    wavelet transform. Explicitly:
    \begin{align}
        [\oper{W}_f v](x,y)=&\scalar{\uir{}{\myhbar}(x,y)v}{f}
        =\int_{\Space{R}{}}
        e^{\pi \rmi \myhbar        (2yx'+xy)}\, v(x'+x)\,
        \bar{f}(x')dx'\nonumber \\
        =&\int_{\Space{R}{}}\textstyle e^{2\pi \rmi \myhbar yx''}\,
        v(x''+\frac{1}{2}x)\,\bar{f}(x''-\frac{1}{2}x)\, dx''.
      \label{eq:Fourier-Wigner}
    \end{align}
    The last expression is known as \emph{Fourier--Wigner transform}%
    \index{Fourier--Wigner transform}%
    \index{transform!Fourier--Wigner}~\citelist{\cite{Folland89}*{\S~1.4}
      \cite{deGosson11a}*{\S~9.2}}.

    For the representation~\eqref{eq:H1-rep-plane-compl} and the
    functional produced by pairing with the Gaussian\index{Gaussian}
    \(\phi(z)=e^{-\pi\myhbar\modulus{z}^2/2}\) the covariant transform
    \(\oper{W}_\phi\) is:
    \begin{align}
        [\oper{W}_\phi f] (z) =&\int_{\Space{C}{}}
        e^{\pi\myhbar  (z\bar{z}'-\bar{z}z')/2} f(z+z')\,
        e^{-\pi\myhbar\modulus{z'}^2/2}\,dz'\wedge d\bar{z}' \nonumber
        \\
      \label{eq:FSB-projection}
                =&\int_{\Space{C}{}}
         f(z'')\,e^{\pi\myhbar z\bar{z}''}\,
        e^{-\pi\myhbar(\modulus{z''}^2+\modulus{z}^2)/2}\,dz''\wedge d\bar{z}'',
    \end{align}
    where \(z''=z+z'\).  This is Fock--Segal--Bargmann (FSB)
    transform%
    \index{Fock--Segal--Bargmann!transform}%
    \index{transform!Fock--Segal--Bargmann}%
    \index{transform!FSB}, it presents the FSB \emph{reproducing kernel}%
    \index{FSB!reproducing kernel}%
    \index{reproducing!kernel!FSB}%
    \index{kernel!reproducing!FSB}
    \(k(z,z'')=e^{\pi\myhbar
      z\bar{z}''}\,
        e^{-\pi\myhbar(\modulus{z''}^2+\modulus{z}^2)/2}\). Note, that
        the second exponent is usually attributed to the weight~\citelist{ \cite{Berezin71a}
      \cite{Folland89}*{\S~1.6} \cite{CoburnIsralowitzLi11a}}. 

    The
    image \(\FSpace{F}{2}(\Space{C}{})\) of~\eqref{eq:FSB-projection}
    is an irreducible invariant subspace of
    \(\FSpace{L}{2}(\Space{C}{})\) with the corresponding orthogonal
    projection:
    \begin{equation}
      \label{eq:FSB-projection-symb}
      P_F: \FSpace{L}{2}(\Space{C}{}) \rightarrow \FSpace{F}{2}(\Space{C}{}),
    \end{equation}
    provided by~\eqref{eq:FSB-projection}. The space
    \(\FSpace{F}{2}(\Space{C}{})\) is characterised by the
    differential equation \((\partial_{\bar{z}}-z) f=0\), which
    follows from Cor.~\ref{co:cauchy-riemann-integ} with the
    distribution \(a(s,x,y)=\delta '_x(s,x,y)-\rmi
    \delta'_y(s,x,y)\). 

  \item For an AHW group%
    \index{abstract Heisenberg--Weyl (AHW) group}%
    \index{group!abstract Heisenberg--Weyl (AHW)}%
    \index{Heisenberg--Weyl! abstract group (AHW)} \(\tilde{G}\),
    \(H\) being its centre and the
    representation~\eqref{eq:AHW-Shroedinger}, any function \(l\) on
    \(G\) produces a pairing appropriate for the induced covariant
    transform, cf.~\eqref{eq:Fourier-Wigner}:
    \begin{equation}
      \label{eq:AHW-cov-trans}
      [\oper{W}_l v](g,\chi)=\int_G
      \chi^k(g')\, f(g+g')\, \bar{l}(g')\,dg'
,
    \end{equation}
    where we integrate over the Haar measure on \(G\).
  \item For \(G=\Space{D}{}\), \(H=M\), the representation
    \(\uir{}{\myh}\)~\eqref{eq:Dynin-represent} is induced by a
    character of the centre \(C\). Thus, any functional can be used
    for the induced covariant transform to \(C\backslash \Space{D}{}\):
    \begin{align*}
        \lefteqn{[\oper{W}_f w](t,u,v,s,x,y) =
          \int_{\Space{H}{}}e^{\pi\rmi\myh(2s't+st+\frac{1}{4}(x'y-xy')t
            +(2x'+x)u +(2y'+y)v)}}\qquad&  \\
        & {} \times \textstyle
        w(s+s'+\frac{1}{2}(x'y-xy'),x+x',y+y')\, \bar{f}(s',x',y')\,ds'\,dx',dy' \\
        =& \int_{\Space{H}{}} \textstyle
        w(s'+\frac{1}{2}s+\frac{1}{4}(x''y-xy''),x''+\frac{1}{2}x,y''+\frac{1}{2}y) \\
        & \qquad{} \times
        \textstyle\bar{f}(s''-\frac{1}{2}s-\frac{1}{4}(x''y-xy''),x''-\frac{1}{2}x,y''-\frac{1}{2}y)\\
        &\qquad{}\times e^{2\pi\rmi\myh(s''t +x''u +y''v)}
        \,ds''\,dx''\,dy''.
    \end{align*}
    A similarity with the Fourier--Wigner transform~\eqref{eq:Fourier-Wigner} is explicit.
  \item 
    \label{it:SU-Bergman-proj}
    For \(G=\SU\), \(H=K\) and the induced
    representation \(\uir{}{}\)~\eqref{eq:SU-action}, a pairing with
    the function \(l_0(z)=1-\modulus{z}^2\), has the property
      \begin{displaymath}
        \scalar{\uir{}{}(h)v}{l_0}=e^{2\rmi\phi}\scalar{v}{l_0},\qquad h=
          \begin{pmatrix}
            e^{\rmi\phi}&0\\0& e^{-\rmi\phi}
          \end{pmatrix}\in K.
      \end{displaymath}
      Thus, \(l_0\) can be used for the induced covariant transform:
      \begin{align}
        \label{eq:bergman-integral}
          [\oper{W}_{0} v](w)=&\int_D \frac{w\bar{z}+1}{\bar{w}z+1}\,
          v\left(\frac{z+w }{\bar{w} z+1 }\right)\,
          (1-\modulus{z}^2)\,\frac{dz\wedge
            d\bar{z}}{\left.(1-\modulus{z}^2)\right.^{2}}\nonumber \\
          =& (1-\modulus{w}^2)\int_{D}\frac{v(\zeta)}{(1-\bar{\zeta}w)^2}\,
          \frac{d\zeta\wedge d\bar{\zeta}}{1-\modulus{\zeta}^2},\qquad\text{where
          } \zeta=\frac{z+w}{\bar{w}z+1}.
      \end{align}
      Up to the factor \(\frac{1-\modulus{w}^2}{1-\modulus{\zeta}^2}\) discussed in
      Example~\ref{it:SU-representation}, 
      this is known as the \emph{Bergman integral}%
    \index{Bergman!reproducing kernel}%
    \index{reproducing!kernel!Bergman}%
    \index{kernel!reproducing!Bergman}~\cite{Coburn12a}. The image
    space \(\FSpace{B}{2}(\Space{D}{})\) of \(\oper{W}_{0}\) is
    \(\SU\)-invariant subspace of \(\FSpace{L}{2}(\Space{D}{})\), which
    is called \emph{Bergman space}%
    \index{Bergman!space}%
    \index{space!Bergman}. The orthogonal projection:
    \begin{equation}
      \label{eq:bergman-proj}
      P: \FSpace{L}{2}(\Space{D}{}) \rightarrow
      \FSpace{B}{2}(\Space{D}{}),
    \end{equation}
    presented by the Bergman integral~\eqref{eq:bergman-integral} is
    called   the   Bergman projection. On
    \(\FSpace{B}{2}(\Space{D}{})\) the
    integral~\eqref{eq:bergman-integral} acts as a reproducing
      formula%
    \index{Bergman!space!reproducing formula}%
    \index{space!Bergman!reproducing formula}%
    \index{reproducing!formula!Bergman space}%
    \index{formula!reproducing!Bergman space}, cf.~\eqref{eq:reproducing-integral}.

    The Bergman space is in the kernel of the differential operator
    \(\frac{z}{1-\modulus{z}^2}-\partial_{\bar{z}}\). For an expression
    of this operator in terms of \(\SU\) and
    Cor.~\ref{co:cauchy-riemann-integ} see
    \citelist{\cite{Kisil97c}*{Ex.~3.7(a)} \cite{Lang85}*{\S~IX.5}}.
  \end{enumerate}
\end{example}

\subsection{Berezin Covariant Symbol}
\label{sec:berez-covar-symb}

An important observation~\cite{Kisil98a} is that, for a representation \(\uir{}{}\) of
\(G\) in a vector space \(V\), we have a representation
\begin{equation}
  \label{eq:double-representation}
  \uir[3]{}{}(g_1,g_2): A \mapsto\uir{}{}(g_2)^{-1}A \uir{}{}(g_1),
  \qquad (g_1,g_2) \in G\times G
\end{equation}
of \(G\times G\) on the space \(B(V)\) of bounded linear operators on \(V\).
\begin{defn}\citelist{\cite{Kisil11c}*{\S~4.3}\cite{Kisil12b}}
  For a fixed operator \(F: B(V)\rightarrow U\) the \emph{covariant
    symbol}%
  \index{covariant!symbol}%
  \index{symbol!covariant} \(\tilde{A}(g_1,g_2)\) is the covariant
  transform defined by the representation \(\uir[3]{}{}\) and the
  operator \(F\):
  \begin{equation}
    \label{eq:covariant-symbol}
    \tilde{A}(g_1,g_2)=F(\uir[3]{}{}(g_1,g_2)
    A)=F(\uir{}{}(g_2)^{-1}A\uir{}{}(g_1)),\qquad
    \text{ where } (g_1,g_2)\in G\times G.
  \end{equation}
  We also use the notation \(\tilde{A}(g)\) for \(\tilde{A}(g,g)\).
\end{defn}
Since the covariant symbol is a special case of the covariant
transform, the respective variants for the scalar case and induced
form are applicable as well. The combination of both has the special
name.  For fixed \(f\in V\) and \(l\in V^*\), the \emph{Berezin
  covariant symbol}%
\index{Berezin!covariant symbol|(}%
\index{covariant!Berezin symbol|(}%
\index{symbol!Berezin covariant|(}%
\index{symbol!Wick|see{Berezin covariant symbol}}%
\index{Wick!symbol|see{Berezin covariant symbol}} \(\tilde{A}(x_1,x_2)\)
is the induced covariant transform defined by the representation
\(\uir[3]{}{}\) and the functional \(F(A)= l(Af)\):
\begin{equation}
  \label{eq:Berezin-symbol}
  \tilde{A}(x_1,x_2)=F(\uir[3]{}{}(\map{s}(x_1),\map{s}(x_2))
  A)=l(\uir{}{}(\uir{}{}(\map{s}(x_2)^{-1})A\map{s}(x_1)) f),
\end{equation}
where  \(x_1\), \(x_2\in H\backslash G\). Again, we denote
\(\tilde{A}(x)=\tilde{A}(x,x)\). 

As before, this definition is most useful if \(f\) and \(l\) are
eigenvectors for all transformations \(\uir{}{}(h)\), \(h\in H\). An
important particular case of the construction is a
unitary representation in a Hilbert space \(V\) and the
functional \(l\in V^*\) be a pairing with \(f\in V\)~\cite{Berezin72}*{\S~1.2}:
\begin{equation}
  \label{eq:Wick-symb-original}
  \tilde{A}(x,y)=\scalar{\uir{}{}(\map{s}(y))^{-1}A\uir{}{}(\map{s}(x))f}{f}=
  \scalar{A\uir{}{}(\map{s}(x))f}{\uir{}{}(\map{s}(y))f}=\scalar{Af_x}{f_y},
\end{equation}
where \(f_x=\uir{}{}(\map{s}(x))f\), \(f_y=\uir{}{}(\map{s}(y))f\).
\begin{example}
  \label{ex:wick-symbol}
  There is a large variety of possibilities (even for a fixed group
  \(G\)) provided by a selection of various subgroups \(H\),
  representations \(\uir{}{}\) and fiducial functionals \(F\). We will
  illustrate this for the Heisenberg group. Note that, our list is
  based on the most popular options and is far from being exhausting.
  For other groups, the number of possibilities is not smaller.
  \begin{enumerate} 
  \item For the Heisenberg group, to make a structure of the listed
    options we introduce a subdivision.
    \begin{enumerate}
    \item \label{it:Berezin-Schwartz-kernel}
      For \(G=\Space{H}{}\) and the
      representation~\eqref{eq:Schroedinger-alt}, take
      \(f(y)=l(y)=\delta(y)\)---the Dirac delta function.  For the
      subgroup \(H_x\)~\eqref{eq:H-x-subgroup} and the homogeneous
      space \(\Space{R}{}=H_x\backslash \Space{H}{}\)
      representation~\eqref{eq:Schroedinger-alt} acts on
      \([\uir{}{}(\map{s}(0,-x,0)
      \delta](x')=\delta_x(x')=\delta(x'-x)\). Consider a smoothing
      operator%
      \index{smoothing!operator}%
      \index{operator!smoothing} \(A: \FSpace[\prime]{S}{}\rightarrow
      \FSpace{S}{}\), where \(\FSpace{S}{}\) is the \emph{Schwartz
        space}%
      \index{Schwartz!space}%
      \index{space!Schwartz} of smooth rapidly decreasing functions on
      the real line and \(\FSpace[\prime]{S}{}\) is its dual---the
      space of tempered distributions. Then the Berezin covariant
      symbol is:
      \begin{equation}
        \label{eq:Berezin-Schwartz}
        \tilde{A}(x_1,x_2)=\scalar{A\delta_{x_1}}{\delta_{x_2}},
      \end{equation}
      which will be related to the Schwartz kernel%
      \index{Schwartz!kernel}%
      \index{kernel!Schwartz} below.
 
      The reader may notice that our usage of the Heisenberg group
      looks excessive in this case: shifts on the real line are
      completely sufficient. Thus, we are moving to the next case.
    \item\label{it:Berezin-KN-symbol} Again consider \(G=\Space{H}{}\)
      this time with the subgroup \(H_x\)~\eqref{eq:H-x-subgroup} and
      the analogous subgroup \(H_y=\{(s,0,x)\in\Space{H}{}\}\).
      Accordingly, for the representation~\eqref{eq:Schroedinger-alt}
      we take \(f(x)\equiv 1\) and \(l(x)=\delta(x)\) both being
      tempered distributions from \(\FSpace[\prime]{S}{}\).  Then,
      \([\uir{}{\myhbar}(0,0,y)f] (x')= e^{\pi\rmi\myhbar 2x'y}
      f(x')\) and \([\uir{}{\myhbar}(0,-x,0)\delta ](x')=
      \delta(x'-x)\). Since, the Fourier transform of
      \(\uir{}{\myhbar}(0,0,y)f\) is the delta function \(\delta_y\),
      for the PDO \(A_{\text{KN}}\)~\eqref{eq:PDO-KN-defn} with a
      smooth symbol \(a\), the Berezin symbol
      \begin{displaymath}
        \tilde{A}_{\text{KN}}(y,x)=a(x,y),
      \end{displaymath}
      is its Kohn--Nirenberg symbol%
      \index{Kohn--Nirenberg symbol}%
      \index{symbol!Kohn--Nirenberg} \(a\).
    \item\label{it:Berezin-Toeplitz} For \(G=\Space{H}{}\), \(H=Z\),
      the representation~\eqref{eq:H1-rep-plane-compl} and the both
      \(l\) and \(f\) be the Gaussian
      \(\phi(z)=e^{-\pi\myhbar\modulus{z}/2}\), the
      transformation~\eqref{eq:Wick-symb-original} is the Wick (or
      Berezin) symbol of an operator \(A\) \citelist{\cite{Berezin71a}
        \cite{Howe80b}
        \cite{Folland89}*{\S~2.7}\cite{CoburnIsralowitzLi11a}}.  The
      simplest calculation of the covariant symbol can be performed
      for the \emph{Toeplitz operator}%
      \index{Toeplitz!operator!FSB space}%
      \index{operator!Toeplitz!FSB space}%
      \index{FSB!operator Toeplitz} \(T_a=P_F a P_F\), with
      \(a\in\FSpace{L}{\infty}(\Space{C}{})\) and
      \(P_F\)~\eqref{eq:FSB-projection-symb}. For the Gaussian
      \(\phi\) and \(\phi_z=\uir{F}{\myhbar}(\map{s}(z))\phi\) we
      found:
     \begin{align}
         \tilde{T}_a(w,z)&=\scalar{T_a \phi_w}{\phi_z}=\scalar{P_F a
           \phi_w}{\phi_z}=\scalar{ a
           \phi_w}{P_F^*\phi_z}=\scalar{a \phi_w}{\phi_z}\nonumber \\
         &=\int_{\Space{C}{}} a(z')\,
         e^{-\pi\myhbar(\bar{w}z'+\modulus{w}^2/2+\modulus{z'}^2/2)}
         e^{-\pi\myhbar(  {z}\bar{z}'+\modulus{z}^2/2+\modulus{z'}^2/2)}\,
         dz'\wedge d\bar{z}' \nonumber \\
       \label{eq:toeplitz-berezin-fock}
         &=e^{-\pi\myhbar(\modulus{w}^2/2+\modulus{z}^2/2)}\int_{\Space{C}{}} a(z')\,
         e^{-\pi\myhbar(\bar{w}z'+{z}\bar{z}'+\modulus{z'}^2)}
         \, dz'\wedge d\bar{z}'.
     \end{align}
     Clearly, \(\tilde{T}_a(z,z)\) is not much different from the FSB
     transform~\eqref{eq:FSB-projection} of \(a\).

     It is worth to notice, that the unitary equivalent model on the
     real line appears if both \(f\) and \(l\) are the Gaussians
     \(e^{-\pi x^2/2}\) on the real line. The respective contravariant
     symbol translates to the language of quantum mechanics as the
     transition amplitude of a quantum mechanical%
     \index{quantum!mechanics}%
     \index{mechanics quantum} observable (in the Schr\"odinger model)
     between states with minimal uncertainty.

     Another class of operators with a useful Berezin calculus are
     \emph{composition operators}%
     \index{composition operator}%
     \index{operator!composition}~\cite{Coburn12a}, i.e. an operator
     \(C_\phi: f \mapsto f\circ \phi\) for a fixed map \(\phi: X
     \rightarrow X\) of the domain to itself.

   \item\label{it:Berezin-localisation} There is another approach for
     \(G=\Space{H}{}\) and the
     representation~\eqref{eq:Schroedinger-alt}.  We take an
     (operator-valued) fiducial operator \(F:
     B(\FSpace{L}{2}(\Space{R}{})) \rightarrow
     B_s(\FSpace{L}{2}(\Space{R}{}))\), where
     \(B_s(\FSpace{L}{2}(\Space{R}{}))\) is the space of bounded
     shift-invariant operators on \(\FSpace{L}{2}(\Space{R}{})\).
     \(F\) is defined by:
    \begin{equation}
      \label{eq:localisation}
      F: A\mapsto A_0,\qquad \text{such that}\quad
      \lim_{\delta \rightarrow 0}|||M_\delta A M_\delta -A_0|||=0,
    \end{equation}
    where \(M_\delta\) is an operator of multiplication by the
    indicator function of \(\delta\)-neighbourhood of the origin and
    \(|||\cdot |||\) denotes the essential norm (modulo compact
    operators). The limit exists for \emph{operators of local type}%
    \index{operator!local type, of}~\cite{Simonenko65a}.
    
    In particular, for the operator \(M_f\) of multiplication by a
    function \(f(x)\) we have \(F M_f=f(0)I\). Therefore, for the
    representation \(\uir{}{\myhbar}\)~\eqref{eq:Schroedinger-alt},
    we obtain the eigenfunction property \(F
    \uir{}{\myhbar}(s,0,y)=e^{2\pi\rmi\myhbar s}I\) for all 
    \((s,0,y)\in H_x\)~\eqref{eq:H-x-subgroup}. Thus, we can use the
    induced form of the covariant symbol~\eqref{eq:covariant-symbol}
    only for values \(\uir{}{\myhbar}(0,x,0)\), where
    \(x\in\Space{R}{}=H_x\backslash
    \Space{H}{}\)~\eqref{eq:Schroedinger-alt}---they are shifts on 
    the real line. Thus the localisation map~\eqref{eq:localisation}
    defines the covariant transform 
    \begin{displaymath}
      A_x=F(\uir{}{\myhbar}(0,x,0) A
    \uir{}{\myhbar}(0,-x,0)),
    \end{displaymath}
    which is the local representative of the operator \(A\) at a point
    \(x\)~\cites{Simonenko65a,Kisil12b}.

    Another important example of operators of local type are
    SIOs\index{SIO}---con\-vo\-lu\-tions on \(\Space{R}{}\) with
    singular kernels---moreover, \(F(S)=S\) for any SIO \(S\) with
    a homogeneous kernel. This recovers Simonenko's localisation
    technique for the calculus of operators generated by SIOs and
    operators of
    multiplications~\cites{Simonenko65a,Simonenko65b,Kisil12b}.
    \end{enumerate}

  \item For an AHW group%
    \index{abstract Heisenberg--Weyl (AHW) group}%
    \index{group!abstract Heisenberg--Weyl (AHW)}%
    \index{Heisenberg--Weyl! abstract group (AHW)} \(\tilde{G}{}\), we
    can essentially repeat all approaches for \(\Space{H}{}\)
    described above. For example, we provide an analogue
    of~\ref{ex:wick-symbol}\ref{it:Berezin-KN-symbol}. For
    \(\tilde{G}\) generated by a commutative group \(G\) consider
    subgroups \(H_G \)~\eqref{eq:AHW-G-subgroup} and the similar
    subgroup \(H_{\hat{G}}=\{(z,g,0)\}\), the respective homogeneous
    spaces are \({G}=H_G\backslash G\) and
    \(\hat{G}=H_{\hat{G}}\backslash G\). Take \(f(g)\equiv 1\) on \(G\) and
    \(l(g)=\delta(g)\). For the
    representation~\eqref{eq:AHW-Shroedinger}, we have
    \([\uir{}{1}(1,0,\chi) f](x')=\chi(x')\) and
    \([\uir{}{1}(1,-x,1) \delta ](x')=\delta_x(x')=\delta (x'-x)\). Then, the Berezin
    symbol of the operator \(A\)~\eqref{eq:AHW-PDO-first} is:
    \begin{displaymath}
      \tilde{A}(x,\chi)=\scalar{A \chi}{\delta_x}=a(x,\chi),
    \end{displaymath}
    i.e. the symbol \(a\) entering the
    integral~\eqref{eq:AHW-Shroedinger}. 
    Other variations
    of~\ref{ex:wick-symbol}\ref{it:Berezin-Schwartz-kernel}--\ref{it:Berezin-localisation} 
    can be obtained in similar ways.

  \item\label{it:Berezin-Dynin-localisation} For \(G=\Space{D}{}\),
    \(H=M\) and the representation~\eqref{eq:Dynin-represent}, we can
    use the localisation approach
    from~\ref{ex:wick-symbol}\ref{it:Berezin-localisation}. For a
    localisation functional \(F\) at the origin of \(\Space{H}{}\)
    similar to~\eqref{eq:localisation}, we calculate
    \(F(\uir{}{\myh}(z,t,u,v,0,0,0))=e^{2\pi\rmi\myh z}I\) for the
    representation~\eqref{eq:Dynin-represent} and \((z,t,u,v,0,0,0)\)
    in the subgroup \(M\)~\eqref{eq:M-subgroup-Dynin}. Thus, it is
    sufficient to perform the covariant
    transform~\eqref{eq:covariant-symbol} for
    \([\uir{}{\myh}(0,0,0,0,s,x,y)\), which are shifts on
    \(\Space{H}{}\). In this way we recovered the calculus of SIO on
    the Heisenberg group initiated in~\cites{Dynin75,Dynin76}, see
    also~\cite{Kisil93b,Kisil94a,Kisil96e}. This can be extended to
    more general nilpotent Lie groups. To this end we need to consider
    a suitable group of dilations, which acts by automorphisms of the
    nilpotent group~\cite{FollStein82}*{\S~1.A}.  Such a covariant
    calculus was recently considered in~\cite{Kisil12b}. This can be
    compared with the standard wavelet technique extended from the
    \(ax+b\) group to the semidirect product of the Heisenberg group
    and the one-dimensional group of its automorphisms acting by
    dilations~\cite{Ishi10a}.
   \item 
     \label{it:SU-toeplitz}
     For \(G=\SU\), \(H=K\) and the
     representation~\eqref{eq:SU-action} we can follow the suit of
     \ref{ex:wick-symbol}\ref{it:Berezin-Toeplitz} by setting
     \(f(z)=l(z)= 1-\modulus{z}^2\). The Berezin
     symbol~\eqref{eq:Wick-symb-original} is
     well-known~\cite{Berezin75b} and very important in the theory of
     operators~\citelist{\cite{Coburn12a}
       \cite{Vasilevski08a}*{\S~A.3} \cite{Nikolski02a}*{\S~B.4.1.8}}.
     Similarly to the Heisenberg group, the simplest calculation of
     the covariant symbol appear for the \emph{Toeplitz operator}%
     \index{Toeplitz!operator}%
     \index{operator!Toeplitz} \(T_a=Pa P\), where
     \(a(z)\in\FSpace{L}{\infty}(\Space{D}{})\) and \(P\) is the Bergman
     projection~\eqref{eq:bergman-proj}. Using expressions from
     Examples~\ref{it:SU-hs-disk} and~\ref{it:SU-representation}, for
     the \(l_0(\zeta)= 1-\modulus{z}^2\) we calculate:
     \begin{displaymath}
      l_w(\zeta)=[\uir{}{}(\map{s}(z))l_0](\zeta)
      =\frac{(1-\modulus{w}^2)(1-\modulus{\zeta}^2)}{(1+\bar{w}\zeta)^2}. 
     \end{displaymath}
     Then:
     \begin{align}
         \tilde{T}_a(w,z)&=\scalar{T_a l_w}{l_z}=\scalar{P a
           l_w}{l_z}=\scalar{ a l_w}{P^* l_z}=\scalar{a l_w}{l_z}\nonumber \\
         &=\int_{\Space{D}{}} a(\zeta)
         \frac{(1-\modulus{w}^2)(1-\modulus{\zeta}^2)}{(1+\bar{w}\zeta)^2}
         \overline{\left(\frac{(1-\modulus{z}^2)(1-\modulus{\zeta}^2)}{(1+\bar{z}\zeta)^2}\right)}\,
         \frac{d\zeta\wedge
          d\bar{\zeta}}{\left.(1-\modulus{\zeta}^2)\right.^2}\nonumber
        \\
       \label{eq:toeplitz-berezin-bergman}
         &=(1-\modulus{w}^2)(1-\modulus{z}^2)\int_{\Space{D}{}} 
         \frac{a(\zeta)}{(1+\bar{w}\zeta)^2(1+{z}\bar{\zeta})^2}\,
         d\zeta\wedge
          d\bar{\zeta}
.
     \end{align}
     Another opportunity to investigate operators on the Bergman space
     is the localisation technique similar
     to~\ref{ex:wick-symbol}\ref{it:Berezin-localisation}. The
     localisation can be combined with the Berezin calculus~\cite{Vasilevski08a}.
  \end{enumerate}
\end{example}%
\index{Berezin!covariant symbol|)}%
\index{covariant!Berezin symbol|)}%
\index{symbol!Berezin covariant|)}%

\subsection{Calculus of Covariant Symbols}
\label{sec:calc-covar-symb}

If a functional \(F\) and a representation \(\uir{}{}\) are both
linear, then the resulting covariant transform
\(\oper{W}_F\)~\eqref{eq:covariant-trans} is a linear map. If
\(\oper{W}_F\) is injective, e.g. due to irreducibility of
\(\uir{}{}\), then \(\oper{W}_F\) transports a norm \(\norm{\cdot}\)
existing on \(V\) to a norm \(\norm[F]{\cdot}\) on the image space
\(\oper{W}_F V\) by the simple rule~\cite{Kisil12d}:
\begin{equation}
  \label{eq:transported-norm}
  \norm[F]{u}:=\norm{v}, \qquad
  \text{ where the unique }
  v\in V \text{ is defined by } u=\oper{W}_F v.
\end{equation}
By the very definition, \(\oper{W}_F\) is an isometry
\((V,\norm{\cdot})\rightarrow (\oper{W}_F V, \norm[F]{\cdot})\).
Moreover, if the representation \(\uir{}{}\) acts on
\((V,\norm{\cdot})\) by isometries then \(\norm[F]{\cdot}\) is right
invariant due to Lem.~\ref{le:covariant-intertwining}.

In most cases, the transported norm can be naturally expressed in the
original terms for \(G\). For example, for a square integrable modulo
a subgroup \(H\) representation%
\index{square integrable!representation}%
\index{representation!square integrable} \(\uir{}{}\) and an
admissible mother wavelet%
\index{wavelet!mother!admissible}%
\index{mother wavelet!admissible}%
\index{admissible mother wavelet} \(f\in V\) the transported
by~\eqref{eq:wavelet-transform} norm coincides with the
\(\FSpace{L}{2}\)-norm on \(X=H\backslash G\). Explicitly, for
\(v_{1,2}\in V\) and
\(\tilde{v}_{1,2}(x)=\scalar[V]{v_{1,2}}{\uir{}{}(\map{s}(x))
  f}\)~\cite{AliAntGaz00}*{Ch.~8}: 
\begin{equation}
  \label{eq:sq-int-isometry}
  \scalar[V]{v_1}{v_2}=\scalar[\oper{W}]{\tilde{v}_1}{\tilde{v}_2},\quad
  \text{where } \scalar[\oper{W}]{\tilde{v}_1}{\tilde{v}_2}=\int_X
  \tilde{v}_1(x)\, \overline{\tilde{v}_2(x)}\,dx.
\end{equation}
Another example of a transported norm is the norm on the Hardy
space in the half-plane~\cite{Kisil12d}.

The particular case of the above transportation is provided by the
Berezin transform. For an operator \(A\) on a normed space \(V\), the
norm of \(A\) has the standard definition:
\(\norm{A}=\sup_{\norm{v}\leq 1} \norm{Av}\). For an isometric
representation \(\uir{}{}\) of \(G\) on \(V\) and \(\norm{f}\leq 1\)
and \(\norm{l}\leq 1\), the associated Berezin transform
\(\tilde{A}(x,y)\)~\eqref{eq:Wick-symb-original} is a function on \(X
\times X\) bounded by \(\norm{A}\). The opposite
statement---boundedness of \(\tilde{A}(x,y)\) implies boundedness of
\(A\)---is a variation of the \emph{reproducing kernel thesis (RKT)}%
\index{reproducing!kernel!thesis (RKT)}%
\index{thesis!reproducing kernel (RKT)}%
\index{RKT|see{reproducing kernel
    thesis}}~\cite{Nikolski02a}*{\S~B.4.1.8}. Another related topic is
a connection of compactness of \(A\) and vanishing of \(A(x,y)\)
``near to the boundary''~\cites{Coburn12a,CoburnIsralowitzLi11a}. We
return to RKT in Subsection~\ref{sec:bound-comp}.

The isometric property~\eqref{eq:sq-int-isometry} allows us to
follow~\cite{Berezin72}*{\S~1.2} and deduce composition rule for Berezin
covariant symbols~\cite{Kisil98a}*{Prop.~3.2}:
\begin{align}
      \widetilde{AB}(x,y)
      &=\scalar[V]{AB f_x}{f_y}=\scalar[V]{B f_x}{A^* f_y}\nonumber \\
      &=\int_X \widetilde{B f_x}(z)\,\overline{\widetilde{A^* f_y}(z)}\,dz
      =\int_X \scalar{B f_x}{f_z} \scalar{f_z}{A^* f_y} dz\nonumber \\
      &=\int_X \scalar{B f_x}{f_z}\scalar{Af_z}{ f_y} dz
      =\int_X \tilde{B}(x,z)\,\tilde{A}(z,y)\,dz.
  \label{eq:Berezin-cov-comp}
\end{align}
One can note, that covariant symbols behaves (up to the order of \(A\)
and \(B\) in the last integral) like integral kernels and this is not
a simple coincidence, see below. Our formula is more
straightforward than the original~\cite{Berezin72}*{\S~1.2} since we
do not need a normalization.

\begin{example}
  \label{ex:covar-compos}
  \begin{enumerate}
  \item For the Heisenberg group in the setup
    of~\ref{ex:wick-symbol}\ref{it:Berezin-Schwartz-kernel}, the fiducial
    functional of pairing with \(f=\delta\) produces the identity
    operator in
    \(\tilde{v}(x)=\scalar{v}{f_x}\)~\eqref{eq:wavelet-transform}.
    Since we have the isometry~\eqref{eq:sq-int-isometry} in the
    trivial way, the composition rule~\eqref{eq:Berezin-cov-comp}
    follows.  Keeping in mind that the Berezin covariant
    transform~\eqref{eq:Berezin-Schwartz} is the Schwartz kernel with
    reversed arguments, we obtained the well-known integral formula for
    the composition of  Schwartz kernels.

    In the setup of~\ref{ex:wick-symbol}\ref{it:Berezin-Toeplitz}, the
    covariant transform turns to be a reproducing
    formula~\eqref{eq:FSB-projection} on the FSB space, thus, is an
    isometry.  The specialisation of the composition
    rule~\eqref{eq:Berezin-cov-comp} for the Toeplitz operators in the
    FSB space can be found in many works starting
    from~\cite{Berezin71a}.
  \item For an AHW group%
    \index{abstract Heisenberg--Weyl (AHW) group}%
    \index{group!abstract Heisenberg--Weyl (AHW)}%
    \index{Heisenberg--Weyl! abstract group (AHW)} \(\tilde{G}{}\), we
    can essentially repeat all approaches (Schwartz kernel, PDO-type,
    Toeplitz operators and localisation techniques) which are in use
    for the Heisenberg group with respective norms and compositions
    formulae. 
  \item For the Dynin group \(\Space{D}{}\) and the
    representation~\eqref{eq:Dynin-represent}, we recall the localisation
    context from~\ref{ex:wick-symbol}\ref{it:Berezin-localisation}
    and~\ref{it:Berezin-Dynin-localisation}. Let \(P_{0,\delta}\) be the
    projection provided by multiplication with the characteristic
    function of \(\delta\)-neigh\-bour\-hood of \(0\in\Space{H}{}\). Then,
    the representation~\eqref{eq:Dynin-represent} produces similar
    projections \(P_{g,\delta}\) for an arbitrary
    \(g\in\Space{H}{}\). For an operator \(A\) on
    \(\FSpace{L}{2}(\Space{H}{})\) we can build a Berezin covariant
    symbol \(\tilde{A}_\delta(g_1,g_2)=P_{g_2,\delta} A
    P_{g_1,\delta}\). If \(A\) is an operator of local type%
    \index{operator!local type, of}~\cite{Simonenko65a}, then
    \(\tilde{A}_\delta(g_1,g_2)\) is a compact for all
    \(\delta<\modulus{g_1-g_2}\). Thus, modulo compact operators the symbol
    \(\tilde{A}(g_1,g_2)=\lim_{\delta\rightarrow
      0}\tilde{A}_\delta(g_1,g_2)\) vanishes outside of the
    diagonal. Therefore, the covariant symbol becomes a field of local
    representatives \(\tilde{A}(g)=\tilde{A}(g,g)\),
    \(g\in\Space{H}{}\)~\cites{Dynin76,Kisil93b}. The
    isometry~\eqref{eq:transported-norm} becomes \(|||A|||=\sup_g
    \norm{\tilde{A}(g)}\). The composition
    rule~\eqref{eq:Berezin-cov-comp} reduces to point-wise
    multiplication of local representatives: \(\widetilde
    {AB}(g)=\tilde{A}(g)\tilde{B}(g)\). 
  \item \label{it:toeplitz-covar-compos} For \(G=\SU\), \(H=K\) and
    the representation~\eqref{eq:SU-action} we also have the
    reproducing formula~\eqref{eq:bergman-integral} on the Bergman
    space. The respective composition formula for Toeplitz operators
    is well-known \citelist{\cite{Berezin75b}*{\S~4.2}
      \cite{Coburn12a} \cite{Vasilevski08a}*{\S~A.3}}.
  \end{enumerate}
\end{example}

\section{Relative Convolutions}
\label{se:relative}

\subsection{Integrated Representations and Contravariant Symbols}
\label{ss:notation}


Let \(G\) be a locally compact group, a left-invariant (Haar) measure on \(G\) is
denoted by \(dg\).  Let \(\uir{}{}\) be a representation of the group
\(G\) in a vector space \(V\). The representation can be extended to a
function \(k\) on \(G\) though integration
\begin{equation}
  \label{eq:integrated-rep}
  \uir{}{}(k)=\int_{G} k(g)\,\uir{}{}(g)\,dg.
\end{equation}
In the simplest case \(k\) has scalar values, however, the same
formula is meaningful for functions with values in operators on the representation
space \(V\). 

The integral~\eqref{eq:integrated-rep} can be defined in a weak sense
for various combinations of functions and representations. One of the
natural setups is a bounded (e.g. unitary) representation \(\uir{}{}\)
and a summable function \(k\). In this case we obtain a
homomorphism of the convolution algebra \(\FSpace{L}{1}(G,dg)\) to an
algebra of bounded operators on \(V\):\index{convolution}
\begin{displaymath}
  \uir{}{}(k_1)\uir{}{}(k_2)=\uir{}{}(k_1*k_2), \qquad\text{ where }\quad
  [k_1*k_2](g)=\int_G k_1(g_1)\,k_2(g_1^{-1}g)\,dg_1.
\end{displaymath}

For a representation \(\uir{}{}\) induced from a subgroup
\(H\), all operators \(\uir{}{}(h)\), \(h\in H\) act
in~\eqref{eq:def-ind} locally. That becomes especially trivial if
\(\uir{}{}(h)\) are scalars. Thus, for induced representations, we are
mainly interested in the ``complement'' \(H\backslash G\) in the
expression~\eqref{eq:integrated-rep}.  For a continuous section
\(\map{s}: H\backslash G \rightarrow G\), we
rewrite~\eqref{eq:integrated-rep} to become an operator of a
\emph{relative convolution}%
\index{relative convolution}%
\index{convolution!relative}~\cite{Kisil94e}:
\begin{equation}
  \label{eq:relative-conv}
  \uir{}{}(k)=\int_{X} k(x)\,\uir{}{}(\map{s}(x))\,dx,
\end{equation}
with a kernel \(k\) defined on \(X=H\backslash G\) with a
(quasi-)invariant measure \(dx\)~\eqref{eq:hs-measure}. Again, the
most natural domain of this definition is a bounded representation
\(\uir{}{}\) and a summable \(k\) from \(\FSpace{L}{1}(X,dx)\).
Furthermore, in many cases we need to (and can) extend meaning
of~\eqref{eq:relative-conv} for suitable functions and distributions,
e.g. the Dirac delta function and its derivatives.

\begin{example}\upshape
  \begin{enumerate}
  \item We already mentioned that relative convolutions generated by
    the Schr\"odinger representation~\eqref{eq:schroedinger-rep} of the
    Heisenberg group are PDO~\eqref{eq:pdo-def}. In this case
    \(G=\Space{H}{}\) and \(H=\{(s,0,0)\such
    s\in\Space{R}{}\}\)---the centre of \(G=\Space{H}{}\).
    It is the original inspiration for this approach~\cites{Howe80b,Folland89,Kisil94e}.
  \item For the AHW group%
    \index{abstract Heisenberg--Weyl (AHW) group}%
    \index{group!abstract Heisenberg--Weyl (AHW)}%
    \index{Heisenberg--Weyl! abstract group (AHW)} \(\tilde{G}\) and
    its representation~\eqref{eq:AHW-Shroedinger} with \(k=1\), take a
    function \(\sigma (g,\chi)\) on \(X=\Space{T}{}\backslash \tilde{G}=G\times \hat{G}\) and calculate,
    cf.~\cite{Folland89}*{(2.32)}:
    \begin{align}
        [\uir{}{}(a) f](g')=&\int_G\int_{\hat{G}} \sigma (g,\chi)\,
        \chi(g')\, f(g+g')\,dg\,d\chi \nonumber  \\
        =&\int_G \hat{\sigma}_2 (g,g')\,
        f(g+g')\,dg \nonumber \\
        =&\int_G \int_{\hat{G}} \hat{\sigma} (\xi,g')\, \bar{\xi}(g)\,
        f(g+g')\,dg\,d\xi \nonumber \\
        =&\int_G \int_{\hat{G}} \hat{\sigma} (\xi,g')\,
        \bar{\xi}(g''-g')\,
        f(g'')\,dg''\,d\xi \nonumber \\
        =&\int_{\hat{G}} \xi(g') \hat{\sigma} (\xi,g')\, \int_G
        f(g'')\,\bar{\xi}(g'') \,dg''\,d\xi,
       \label{eq:AHW-PDO-second}
    \end{align}
    here \(g''=g+g'\), \(\hat{\sigma}_2\) is a function on \(G\times
    G\), which is the Fourier transform of \(\sigma\) in second
    variable. The last expression~\eqref{eq:AHW-PDO-second} coincides with 
    Kohn--Nirenberg type PDO~\eqref{eq:AHW-PDO-first} for
    \(a(g,\xi)=\hat{\sigma}(\xi,g)\), cf.~\amscite{RuzhanskyTurunen10a}*{Part~II}.
  \item 
    \label{it:dynin-source} 
    For the Dynin group \(\Space{D}{}\), the unitary representation
    \(\uir{}{}\)~\eqref{eq:Dynin-represent} on
    \(\FSpace{L}{2}(\Space{H}{})\) is obtained from its infinitesimal
    action~\eqref{eq:h-lie-algebra} and~\eqref{eq:dynin-mult-rep}. The
    integrated representation~\eqref{eq:integrated-rep} was considered
    in~\cite{Dynin75} as a generalisation of the Weyl quantization
    from \(\Space{H}{}\) to the group \(\Space{D}{}\).

    If function \(k\) has the structure
    \(k(z,t,u,v,s,x,y)=\delta(z,t,u,v)k_1(s,x,y)\), where \(\delta\)
    is the Dirac delta function, then \(\uir{}{}(k)\) is a convolution
    on the Heisenberg group with the kernel \(k_1\). On the other
    hand, if
    \begin{displaymath}
      k(z,t,u,v,s,x,y)=\delta(z) k_2(t,u,v) \delta (s,x,y),
    \end{displaymath}
    then \(\uir{}{}(k)\) is an operator of multiplication by
    \(\hat{k}_2(s,x,y)\)---the (Euclidean) Fourier transform
    \((t,u,v)\rightarrow (s,x,y)\) of \(k_2\). Thus, the integrated
    representation~\eqref{eq:integrated-rep} in this case belongs to
    the algebra of operators generated by convolutions on the
    Heisenberg group and multiplications by functions, which were
    investigated, for example,
    in~\cites{Dynin75,Dynin76,Kisil93b,Kisil12b}. For a suitable
    choice of symbols, this operators coincide
    with~\eqref{eq:oper-dual-gr} used
    in~\cite{BahouriFermanian-KammererGallagher12}.

    Furthermore, we can observe that in both cases kernels
    depend on the coordinate \(z\) through the delta function. Thus,
    instead of the integrated representation~\eqref{eq:integrated-rep}
    we can use the relative convolutions~\eqref{eq:relative-conv} for
    \(G=\Space{D}{}\) and \(H\) being its centre, cf. the case of the
    Heisenberg group above. 

  \item For \(G=\SU\) and \(H=K\), a
    substitution of~\eqref{eq:SU-maps} into
    representation~\eqref{eq:SU-action} produces the
    relative convolutions:
    \begin{align}
        [\uir{}{}(k) v](z) =& \int_{D}k(w)\,\frac{w\bar{z}+1
        }{\bar{w}z+1}\, v\left(\frac{z+w }{\bar{w} z+1
          }\right)\frac{dw\wedge
          d\bar{w}}{\left.(1-\modulus{w}^2)\right.^2}\nonumber \\
        =& \int_{D}\frac{1-z\bar{\zeta}}{1-\bar{z}\zeta}\,
        k\left(\frac{z-\zeta }{\bar{z} \zeta-1 }\right)v(\zeta)
        \frac{d\zeta\wedge
          d\bar{\zeta}}{\left.(1-\modulus{\zeta}^2)\right.^2},\quad\text{where
        } w=\frac{z-\zeta}{\bar{z}\zeta-1}.
      \label{eq:SU-rel-conv}
    \end{align}
    Interestingly, the last integral can be interpreted as
    \(\tilde{v}(z)=\scalar{v}{\uir{}{}(\map{s}^{-1}(z))\bar{k}}\), which is
    the induced wavelet transform on
    \(K\backslash\SU\)~\citelist{\amscite{Kisil11c}*{\S~5.5}} with the
    mother wavelet \(\bar{k}\).
  \end{enumerate}
\end{example}
The indicated connection of relative convolutions with the induced
wavelet transform is not an exception. It occurs in many other cases
when the representation space \(V\) consists of functions defined on
the homogeneous space \(X=H\backslash G\), e.g. 
the FSB space of analytic functions on \(\Space{C}{n}=Z\backslash
\Space{H}{n}\). In general, the covariant transform, the Berezin
symbol, integrated representations and the contravariant symbol%
\index{Berezin!contravariant symbol}%
\index{contravariant!Berezin symbol}%
\index{symbol!Berezin contravariant}%
\index{symbol!anti-Wick|see{Berezin contravariant symbol}}%
\index{anti-Wick symbol|see{Berezin contravariant symbol}}
(considered below) are closely connected and, sometimes, even
confused. 

\subsection{Twisted Convolutions}
\label{sec:twisted-convolutions}

It is desirable to have an efficient symbolic calculus of relative
convolutions. For exponential Lie groups, a calculus in terms of the
respective Lie algebras was initiated in~\cite{Kisil94e}. However, the
exponential property is rather restrictive, for example, \(\SU\) does
not posses it. Here we provide another algebraic condition, which is
sufficient for relative convolutions to be closed under
multiplication.

In the notations of Section~\ref{sec:induc-repr}, for any \(x_1\),
\(x_2\in X=H\backslash G\) there is the unique \(x\in X\) defined by
the identity
\begin{equation}
  \label{eq:semugroup-def}
  \map{s}(x_1)\,\map{s}(x_2)=h\map{s}(x), \qquad \text{ that is }
  \quad x=\map{p}(\map{s}(x_1)\,\map{s}(x_2))=x_1\cdot\map{s}(x_2),
\end{equation}
where the last expression uses notation~\eqref{eq:g-action}.

The relation~\eqref{eq:semugroup-def} defines a binary operation
\((x_1,x_2)\mapsto x\), which turns \(X\) into a semigroup. It is not
a group unless \(H\) is a normal subgroup of \(G\). One can develop a
separate theory for semigroups from homogeneous spaces, for example,
in~\cite{Ungar09a} they are called \emph{gyrogroups}. However, we
prefer to proceed in terms of the original group \(G\) and its
subgroup \(H\).

For given \(x_2\), \(x\in X\), there is the only
\(x_1=x\cdot(\map{s}(x_2))^{-1}\in X\) satisfying the
first identity in~\eqref{eq:semugroup-def}. We will use the
abbreviation \(x_1=xx_2^{-1}\) for it.

Furthermore, using the transformation rule~\eqref{eq:hs-measure-transf}
of the measure \(dx\) on \(X\) we calculate:
\begin{displaymath}
  dx_1\,dx_2=\frac{\Delta_H(h(x,x_2))}{\Delta_G(h(x,x_2))}\,dx_2\,dx,\quad
  \text{where} \quad
  h(x,x_2)=\map{s}(x)\map{s}^{-1}(x_2)\map{s}^{-1}(x\cdot \map{s}^{-1}(x_2)).
\end{displaymath}
Here, for simplicity, we write \(\map{s}^{-1}(y)\) instead of the more correct
expression \((\map{s}(y))^{-1}\).

Let two relative convolutions be defined by scalar-valued summable kernels \(k_1\),
\(k_2\in\FSpace{L}{1}(X)\). Then, starting with the Fubini theorem we calculate:
\begin{align}
  \uir{}{}(k_1)  \uir{}{}(k_2)=&
 \int_{X} k_1(x_1)\,\uir{}{}(\map{s}(x_1))\,dx_1  \int_{X}
 k_2(x_2)\,\uir{}{}(\map{s}(x_2))\,dx_2
\nonumber \\
=& 
 \int_{X} \int_{X} k_1(x_1)\,
 k_2(x_2)\,\uir{}{}(\map{s}(x_1)\map{s}(x_2))\,dx_1\,dx_2 \nonumber \\
=& 
 \int_{X} \int_{X} k_1(x_1)\,
 k_2(x_2)\,\uir{}{}(\map{s}(x_1)\map{s}(x_2))\,dx_1\,dx_2\nonumber \\
=& 
 \int_{X} \int_{X} k_1(xx_2^{-1})\,
 k_2(x_2)\,\uir{}{}(h^{-1}(x,x_2)\map{s}(x))\,\frac{\Delta_H(h(x,x_2))}{\Delta_G(h(x,x_2))}\,dx_2\,dx \nonumber \\
=& 
 \int_{X} \int_{X} k_1(xx_2^{-1})\,
 k_2(x_2)\,\uir{}{}(h^{-1}(x,x_2))\,\frac{\Delta_H(h(x,x_2))}{\Delta_G(h(x,x_2))}\,dx_2 
 \,\uir{}{}(\map{s}(x))\,dx \nonumber \\
 =& 
 \int_{X} \int_{X} k_1(xx_2^{-1})\,
 k_2(x_2)\,\uir{}{}(h^{-1}(x,x_2))\,\frac{\Delta_H(h(x,x_2))}{\Delta_G(h(x,x_2))}\,dx_2 
 \,\uir{}{}(\map{s}(x))\,dx \nonumber \\
 =&
 \int_{X}  k(x) \,\uir{}{}(\map{s}(x))\,dx,
 \label{eq:rel-conv-comp}
\end{align}
where
\begin{equation}
  \label{eq:twisted-convolution}
  k(x)=\int_{X} k_1(xx_2^{-1})\,
 k_2(x_2)\,\uir{}{}(h^{-1}(x,x_2))\,\frac{\Delta_H(h(x,x_2))}{\Delta_G(h(x,x_2))}\,dx_2.
\end{equation}
Note that, if the representation \(\uir{}{}\) is induced by a
character \(\chi\) of the subgroup \(H\), then
\(\uir{}{}(h^{-1}(x,x_2))=\chi(h^{-1}(x,x_2))\). Thus, the above
integral is scalar valued.
\begin{defn}
  \label{de:twisted-convolution}
  For two summable functions \(k_1\) and \(k_2\) on \(X=H\backslash
  G\), their \emph{twisted convolution}%
  \index{twisted convolution}%
  \index{convolution!twisted} \(k=k_1 \natural k_2\) is a
  function \(k\), such that the relative convolution
  \(\uir{}{}(k)\)~\eqref{eq:relative-conv} equal to the composition
  \(\uir{}{}(k_1) \uir{}{}(k_2)\) of relative convolutions with the
  kernels \(k_1\) and \(k_2\), i.e.:
  \begin{displaymath}
    \uir{}{}(k_1 \natural k_2)=\uir{}{}(k_1)  \uir{}{}(k_2).
  \end{displaymath}
\end{defn}
Then, the result of
calculations~\eqref{eq:rel-conv-comp} can be encapsulated in the
statement:
\begin{prop}
  The twisted convolutions \(k=k_1 \natural k_2\) of two functions
  \(k_1\) and \(k_2\) on \(H\backslash G\) is presented
  by~\eqref{eq:twisted-convolution}.
\end{prop}
\begin{rem}
  Integrated representations~\eqref{eq:integrated-rep} and
  relative convolutions~\eqref{eq:relative-conv} map
  functions to operators. It is a fashion now to call any such map a
  ``quantization''.  An opposite procedure, e.g. the covariant
  transform, maps an operator to a function---a symbol of the
  operator. This can be called ``dequantization'', respectively. Thus
  our Defn.~\ref{de:twisted-convolution} can be stated in
  quasi-quantum language as follows: quantize kernels to operators,
  compose operators, dequantize the composition to the kernel. In this
  setup the twisted convolution is also known as a \emph{star product}%
  \index{star!product}%
  \index{product!star}%
  \index{star!product|see{twisted convolution}}%
  \index{product!star|see{twisted convolution}}. We refer to~\cite{Aniello09a}
  for further discussion, references and a more explicit formula in the
  case of square integrable representations.%
\index{square integrable!representations}%
\index{representations!square integrable} We also note, that a
  search of a compatible star product for an arbitrary Poisson
  manifold is the topic of \emph{deformation quantization}%
  \index{deformation quantization}%
  \index{quantization!deformation}.
\end{rem}

\begin{example}
  \begin{enumerate}
  \item The Heisenberg group \(\Space{H}{}\) and its centre \(Z\) are
    unimodular, thus \(\Delta_{\Space{H}{}}\equiv 1\) and
    \(\Delta_Z\equiv 1\).  Using maps \(\map{p}\) and \(\map{s}\) from
    Example~\ref{it:heisenberg-map} for \(\Space{R}{2}=Z\backslash
    \Space{H}{}\), we calculate:
    \begin{align*}
        (x,y)(x_2,y_2)^{-1}=&(x-x_2,y-y_2),\\
        h((x,y),(x_2,y_2))=&(\textstyle\frac{1}{2}(x_2y-yx_2),0,0).
    \end{align*}
    Thus, for a representations~\eqref{eq:H1-rep-plane}
    and~\eqref{eq:Schroedinger-alt} induced by a character
    \(\chi_\myhbar(s,0,0)=e^{2\pi\rmi\myhbar s}\), the respective
    twisted convolution is:
    \begin{displaymath}
      (k_1\natural k_2)(x,y)=\int_{\Space{R}{2}} k_1(x-x_2,y-y_2)\,
      k_2(x_2,y_2)\,e^{-\pi\rmi\myhbar (x_2y-yx_2)}\,dx_2\,dy_2.
    \end{displaymath}
    This operation is the key for the whole calculus of PDO as
    explained in~\citelist{\amscite{Howe80b}*{\S~2}
      \amscite{Folland89}*{\S~2.3}}. It is also known as the
    Groenewold--Moyal star product~\cite{Aniello09a}*{\S~6.1}%
  \index{star!product!Groenewold--Moyal}%
  \index{product!star!Groenewold--Moyal}%
  \index{Groenewold--Moyal star product}.
  \item For an AHW group%
\index{abstract Heisenberg--Weyl (AHW) group}%
\index{group!abstract Heisenberg--Weyl (AHW)}%
\index{Heisenberg--Weyl! abstract group (AHW)} \(\tilde{G}\), we calculate
    \(h((g,\chi),(g_2,\chi_2))=\chi_2(g_2-g)\). From unimodularity of
    \(\tilde{G}\), the twisted convolution is:
    \begin{displaymath}
      (k_1\natural k_2)(g,\chi)=\int_G\int_{\hat{G}}
      k_1(g-g_2,\chi-\chi_2)\,k_2(g_2,\chi_2)\,\chi_2(g_2-g)\,dg_2\,d\chi_2.
    \end{displaymath}
    This again looks like a convolution on the Cartesian product
    \(G\times \hat{G}\) with a ``twist''. A composition with the
    Fourier transform on \(G\times \hat{G}\) maps our twisted
    convolution to the star product from~\cite{Wildberger05b}*{\S~3}.
  \item According to Kirillov's theory, any unitary irreducible
    representation of a nilpotent group Lie group is induced by a
    character of the group's centre
    \(C\)~\amscite{Kirillov76}*{\S~15}. Thus, the relative convolution
    for \(C\backslash\Space{D}{}\) is not much different from the
    whole integrated representation.
    
    On the other hand, since the subgroup
    \(M\)~\eqref{eq:M-subgroup-Dynin} is normal then the twisted
    convolution for \(M\backslash \Space{D}{}=\Space{H}{}\) reduces to
    the group convolution on the Heisenberg group. Some interesting
    options are located between the extremes \(C\) and \(M\). For
    example, since the Lie algebra of \(\Space{D}{}\) is generated by
    \(X\), \(Y\), \(T\), it is worth to consider twisted convolution
    associated to the subgroup  
    \begin{displaymath}
      M'=\{(z,0,u,v,0,0,0)\in \Space{D}{} \such (z,u,v)\in\Space{R}{3}\}.
    \end{displaymath}

  \item The group \(G=\SU\) is not unimodular, however
    \(\Delta_{G}(k)\equiv 1\) for all \(k\) in the maximal compact
    subgroup \(K\)~\amscite{Lang85}*{\S~III.1}. The subgroup \(K\) is
    unimodular since it is commutative. Furthermore, using
    maps~\eqref{eq:SU-maps} we calculate: 
    \begin{displaymath}
      z\cdot\map{s}^{-1}(w)=\frac{z-w}{1-z\bar{w}}, \qquad
      h(z,w)=\frac{\modulus{1-\bar{z}w}}{1-\bar{z}w}.
    \end{displaymath}
    Thus, for the representation~\eqref{eq:SU-action} induced from
    \(K\), the respective twisted convolution is:
    \begin{displaymath}
      (k_1\natural k_2)(z)=\int_{D} k_1\left(\frac{z-w}{1-z\bar{w}}\right)\,
      k_2(w)\,\frac{1-\bar{z}w}{1-z\bar{w}}\,\frac{dw\wedge d\bar{w}}{\left.(1-\modulus{w}^2)\right.^2}.
    \end{displaymath}
    This corresponds to the composition of the Berezin contravariant
    symbols%
\index{Berezin!contravariant symbol}%
\index{contravariant!Berezin symbol}%
\index{symbol!Berezin contravariant} (see below) and nicely complements the well-known calculus
    of Berezin's covariant (Wick) symbols considered in
    Example~\ref{it:toeplitz-covar-compos} and
    \citelist{\cite{Berezin75b}*{\S~4.2} \cite{Coburn12a}
      \cite{Vasilevski08a}*{\S~A.3}}.
  \end{enumerate}
\end{example}

\subsection{Contravariant Symbol and Toeplitz Operators}
\label{sec:bound-comp}

There is a notion immediately derived from the integrated
representations: the \emph{contravariant}%
\index{contravariant!transform}%
\index{transform!contravariant} (aka inverse covariant)
transform~\cites{Kisil11c,Kisil10c,Kisil09d,Kisil12d}. For an
integrated representation \(\uir{}{}\)~\eqref{eq:integrated-rep}
(or~\eqref{eq:relative-conv}) and a fixed vector \(w\in V\) the
associated contravariant transform of a function \(k\) on \(G\) (or
\(X=H\backslash G\)) is
\begin{equation}
  \label{eq:contravariant-transform}
  \oper{M}_w^{\uir{}{}} (k)= \uir{}{}(\check{k})w, \qquad \text{ where }
  \quad \check{k}(g)=k(g^{-1}).
\end{equation}
The contravariant transform \(\oper{M}_{w}^{\uir{}{}}\) intertwines the
right regular representation \( R \)
on \( \FSpace{L}{2}(G)\) and \( \uir{}{} \):
\begin{equation}
  \label{eq:inv-transform-intertwine}
  \oper{M}_{w}^{\uir{}{}}\, R(g) = \uir{}{}(g)\, \oper{M}_{w}^{\uir{}{}}.
\end{equation}
Combining with~\eqref{eq:covariant-intertwine}, we see that the
composition \(\oper{M}_w^{\uir{}{}} \circ \oper{W}_v^{\uir{}{}}\) of
the covariant and contravariant transform intertwines \(\uir{}{}\)
with itself.  We can use the Schur's lemma
\citelist{\cite{AliAntGaz00}*{Lem.~4.3.1}
  \cite{Kirillov76}*{Thm.~8.2.1}} to deduce that:
\begin{prop}
  \label{pr:composition}
  For an irreducible \(\uir{}{}\), the composition \(\oper{M}_w^{\uir{}{}} \circ
  \oper{W}_v^{\uir{}{}}\) is a multiple \(\theta(w,v)I\) of the identity
  operator. Moreover, the factor \(\theta(w,v)\) is a sesquilinear form
  of vectors \(w\), \(v\in V\).
\end{prop}
The following interesting consequence is known in slightly different
form for the case of the Heisenberg group~\cite{Folland89}*{(1.47)}.
\begin{cor}
  Assume that the integrated representation \(\uir{}{}\) is faithful
  on the image space \(\oper{W}_{v_2} V\). Then the twisted
  convolution of wavelet transforms is:
  \begin{equation}
    \label{eq:twist-conv-wave-trans}
    \oper{W}_{v_1} u_1 \natural
    \oper{W}_{v_2} u_2= \theta(u_2,v_1)  \oper{W}_{v_2}
    u_1
  \end{equation}
\end{cor}
\begin{proof}
  We note another form \(\uir{}{}(\oper{W}_v u)w=\theta(w,v) u\) of the
  identity \(\oper{M}_w^{\uir{}{}} \circ
  \oper{W}_v^{\uir{}{}}=\theta(w,v)I\). Then:
  \begin{align*}
    \oper{M}_{v_2}(\oper{W}_{v_1} u_1 \natural
    \oper{W}_{v_2} u_2)&= \uir{}{}(\oper{W}_{v_1} u_1 \natural
    \oper{W}_{v_2} u_2) v_2\\
    &= \uir{}{}(\oper{W}_{v_1} u_1) 
    \uir{}{}(\oper{W}_{v_2} u_2) v_2\\
    &= \uir{}{}(\oper{W}_{v_1} u_1) 
    \theta(v_2,v_2) u_2 \\
    &= \theta(v_2,v_2) \uir{}{}(\oper{W}_{v_1} u_1) 
     u_2 \\
    &= \theta(v_2,v_2) \theta(u_2,v_1)  u_1.
  \end{align*}
  Since the representation \(\uir{}{}\) is faithful on the image space
  \(\oper{W}_{v_2} V\),  the obtained
  result implies~\eqref{eq:twist-conv-wave-trans}.
\end{proof}
The following particular case of~\eqref{eq:twist-conv-wave-trans} is
of special interest:
\begin{cor}
  \label{co:reproducing-formula}
  Under the assumptions of the previous Corollary, the image space
  \(\oper{W}_{f} V\) is reproducing kernel%
  \index{reproducing kernel}%
  \index{kernel!reproducing} space with the following
  realisation of a reproducing formula:
  \begin{equation}
    \label{eq:reproducing-formula}
    \tilde{v}= \tilde{f}  \natural \tilde{v}, \qquad
    \text{ for }\quad \tilde{f}=\oper{W}_{f}f \quad\text{ and any }\quad \tilde{w}\in \oper{W}_{f} V.
  \end{equation}
\end{cor}

The contravariant transform is a source of the Berezin's
contravariant symbol%
\index{Berezin!contravariant symbol}%
\index{contravariant!Berezin symbol}%
\index{symbol!Berezin contravariant} as follows. For a pair \(v\in V\), \(f\in V^*\), consider a
rank-one operator \(E_{v,f}:V\rightarrow V\) define by the
expression \(E_{v,f} u=\scalar{u}{f}v\). Then, the representation
\(\uir[3]{}{}\)~\eqref{eq:double-representation} acts as follows:
\begin{equation}
  \label{eq:rank-one-oper}
  \uir[3]{}{}(g_1,g_2) E_{v,f} = E_{v',f'}, \quad \text{ where }\quad
  v'=\uir{}{}(g_2)^{-1}v \quad\text{ and }\quad f'= \uir{*}{}(g_1)f,
\end{equation}
with the last identity natural meaning
\(\scalar{u}{f'}=\scalar{u}{\uir{*}{}(g_1)f}=\scalar{\uir{}{}(g_1)u}{f}\).
Then, the contravariant transform  for the
vector \(w=E_{v,f}\in \FSpace{B}{}(V)\) becomes:
\begin{equation}
  \label{eq:contravariant-symbol}
  \oper{M}_{v,f}(a)=\int_X \int_X a(x_1,x_2)\,
  \uir[3]{}{}(\map{s}(x_1),\map{s}(x_2)) 
  E_{v,f}\,dx_1\,dx_2.  
\end{equation}
Using~\eqref{eq:rank-one-oper}, we can directly write the action of
operator~\eqref{eq:contravariant-symbol} on \(u\in V\):
\begin{equation}
  \label{eq:contravariant-explicit}
  \oper{M}_{v,f}(a) u=\int_X \int_X a(x_1,x_2)\,
  \scalar{\uir{}{}(\map{s}(x_1))u}{f}\,\uir[3]{}{}(\map{s}(x_2)^{-1})v\,dx_1\,dx_2. 
\end{equation}
Here we call  \(a\) the \emph{contravariant symbol}%
\index{contravariant!symbol}%
\index{symbol!contravariant} of the operator \(\oper{M}_{v,f}(a)\).

The contravariant symbol in the sense of
Berezin~\cites{Berezin72,Berezin71a,Berezin75b} appears if \(V\) is a
Hilbert space with an irreducible square integrable representation%
\index{square integrable!representations}%
\index{representations!square integrable}
\(\uir{}{}\). The covariant transform \(\oper{W}_v\) for an admissible
mother wavelet \(v\in V\) identifies \(V\) with its image \(\oper{W}_v
V\). There is the respective reproducing kernel
\(k_y\)~\eqref{eq:reproducing-integral} on \(\oper{W}_v V\).  If the
representation \(\uir[3]{}{}\) is restricted to the diagonal
\(\uir[3]{}{}(g)=\uir[3]{}{}(g,g)\) of \(G\times G\), then the contravariant
transform similar to~\eqref{eq:contravariant-explicit}
is~\cite{Berezin72}*{\S~1.1}:
\begin{align}
    [\oper{M}_{v,v}(a) u](y)&=\int_X a(x)\,
    \scalar{\uir{}{}(\map{s}(x))u}{v}\,\uir[3]{}{}(\map{s}(x)^{-1})v(y)\,dx
     \nonumber \\
    &=\int_X a(x)\,u(x)\, k_y(x)\,dx.
  \label{eq:toeplitz-berezin}
\end{align}
The last expression is the \emph{Toeplitz operator}%
\index{Toeplitz!operator}%
\index{operator!Toeplitz} \(T_a=PaP\) for the projection \(P\) on
\(\oper{W}_v\) defined by the integral in the right-hand side of the
reproducing formula~\eqref{eq:reproducing-integral}.

The explicit formulae connectioning co- and contravariant
symbols%
\index{Berezin!contravariant symbol}%
\index{contravariant!Berezin symbol}%
\index{symbol!Berezin covariant}%
\index{Berezin!covariant symbol}%
\index{covariant!Berezin symbol}%
\index{symbol!Berezin contravariant} are known for a long time~\citelist{\cite{Berezin72}*{(1.12)}
  \cite{Berezin75b}*{(3.13)}}. Within our approach they are
consequences of Prop.~\ref{pr:composition}, since covariant and
contravariant symbols are special cases of covariant and contravariant
transforms.

The original Berezin's papers~\cites{Berezin72,Berezin71a,Berezin75b}
(as well as subsequent developments in the context of abstract
reproducing kernel spaces)%
\index{space!reproducing kernel}%
\index{reproducing!kernel!space}%
\index{kernel!reproducing} do not assume any group structure.  It is
possible to obtain important estimations for norm and compactness in
this abstract setting. However, the fundamental examples---the Bergman
and FSB spaces---considered in those papers are generated by groups,
as we have seen above. In particular, the group structure becomes very
relevant in the study composition operators%
\index{composition operator}%
\index{operator!composition} generated by an automorphism of the
domain~\cite{Coburn12a}. Furthermore, the formula for twisted
convolution~\eqref{eq:twisted-convolution} is also based on the
underlined group structure and is not possible on a generic
reproducing kernel space.

\begin{example}
  \begin{enumerate}
  \item For Toeplitz operators~\eqref{eq:toeplitz-berezin-fock} on the
    Heisenberg group, the contravariant calculus was already
    investigated in~\cite{Berezin71a} and is still an important
    tool~\cites{CoburnIsralowitzLi11a,Coburn12a}. The connections
    between PDO~\eqref{eq:PDO-Weyl-defn} and the Toeplitz
    operators~\eqref{eq:toeplitz-berezin-fock} was fruitfully
    exploited in~\cite{Howe80b}*{\S~4.2} with the following
    observation: ``The Toeplitz operators are to the Bargmann--Fock
    model as the pseudodifferential operators\index{PDO} are to the
    Schr\"odinger model''. We used such a technique to obtain 
    Calder\'on--Vaillancourt--type estimations for relative
    convolutions on exponential nilpotent Lie groups~\cite{Kisil13b}.
  \item For the AHW group \(\tilde{G}\), the connection between
    PDO-type operators~\eqref{eq:AHW-PDO-second}
    \cite{RuzhanskyTurunen10a} and Toeplitz-type operators on
    \(\FSpace{L}{2}(G\times \hat{G})\) shall closely follow the
    Heisenberg group suit. Yet, no work in this direction is known to
    me.
  \item For the Dynin group, I am not aware of any study of
    contravariant calculus and Toeplitz-type operators. However, it is
    natural to expect, that their relation to the PDO-like calculus
    from Example~\ref{it:dynin-source}
    and~\cites{Dynin75,Dynin76,BahouriFermanian-KammererGallagher12,FischerRuzhansky12a}
    shall be similar to the Heisenberg group. However, due to a higher
    level of non-commutativity it may be not as straightforward as for
    an AHW group.
  \item For \(\SU\) and Toeplitz
    operators~\eqref{eq:toeplitz-berezin-bergman}, the Berezin
    contravariant symbols was studied in~\cite{Berezin75b}, with
    numerous fruitful developments,
    cf.~\cites{Coburn12a,Coburn94a,Coburn90,Vasilevski08a}.
  \end{enumerate}
\end{example}

\section{Discussion}
\label{sec:discussion}
The moral of the present overview is that there is no a single formula
perfectly serving all situations. However, covariant and contravariant
transforms provide a general framework which has a rich and  flexible
inventory.  The presented list of different examples prompts further
detailed investigation of this approach in various directions.

\section*{Acknowledgement}
\label{sec:acknowledgement}
I am grateful to Prof.~Gerald B.~Folland for useful comments on the
theory of AHW groups%
\index{abstract Heisenberg--Weyl (AHW) group}%
\index{group!abstract Heisenberg--Weyl (AHW)}%
\index{Heisenberg--Weyl! abstract group (AHW)}, which
pointed to the respective references~\cites{Mackey49a,Mackey70a}.

\appendix
\section{Covariant Transform: A Road Map} 
\label{se:appendix}
 
This paper's presentation was illustrated by numerous detailed
examples. Here we provide a brief outline of notions and main formulae used in
our constructions.

\begin{itemize}
\item Induced representations.
  \begin{itemize}
  \item \(G\) is a locally compact group with a right invariant
    measure \(dg\) and the modular function \(\Delta_G\).
  \item \(H\subset G\) is a subgroup with a right invariant measure
    \(dh\) and the modular function \(\Delta_H\).
  \item \(X=H\backslash G\) is the homogeneous space of right cosets:
    \(g_1\sim g_2\) if \(g_1=hg_2\) for \(h\in H\).
  \item \(\map{p}: G \rightarrow H\backslash G\) is the natural
    projection of an element to its coset.
  \item \(\map{s}: H\backslash G \rightarrow H\) is a section---a right
    inverse of \(\map{p}\): \(\map{p}(\map{s}x)=x\) for all \(x\in
    H\backslash G\). The map \(\map{s}\) is not unique and we often can
    chose it continuous.
  \item \(\map{r}: G \rightarrow H\) is defined from the identity:
    \(g=\map{r}(g)\map{s}(\map{p}(g))\), \(g \in G\).
  \item \(X=H\backslash G\) is a right \(G\)-space with the action: \(g:
    x\mapsto x\cdot g=\map{p}(\map{s}(x)*g)\), \(g\in G\), \(x\in X\).
  \item A representation of \(G\) induced by a character \(\chi\) of the
    subgroup \(H\) is \( [\uir{}{\chi}(g) f](x)=
    \chi(\map{r}(\map{s}(x)*g))\, f(x\cdot g)\).
  \item An equivalent form of the induced representation is by the right
    shift \(R(g): F(g')\mapsto F(g'g)\) on a space of functions with the
    property \(F(hg)=\chi(h)F(g)\), for \(h\in H\) and \(g\in G\).
  \item Two models of induced representations are connected by the lifting
    \([\oper{L}_\chi f](g)=\chi(h) f(\map{p}(g))\) and pulling
    \([\oper{P}F](x)=F(\map{s}(x))\).
  \end{itemize}
\item Covariant transform.
  \begin{itemize}
  \item For a representation \(\uir{}{}\) of \(G\) in \(V\) and an
    operator \(F: V \rightarrow U\), the covariant transform is
    \([\oper{W}_F v](g)=F(\uir{}{}(g)v)\), where \(v\in V\) and \(g\in
    G\).
  \item The induced covariant transform is \( [\oper{W}_F
    v](x)=F(\uir{}{}(\map{s}(x))v)\) for \(v\in V\), \(x\in H\backslash
    G\).
  \item If a mother wavelet \(f\) satisfies to the identity \(\int_{G}
    a(g)\, \uir{*}{}(g) f\,dg=0\) for some distribution \(a(g)\) on \(G\),
    then any wavelet transforms \(\tilde{v}=\scalar{v}{\uir{*}{}(g)f}\)
    satisfies the identity \(D \tilde{v}=0\), where \(D=\int_{G}
    \bar{a}(g)\, \Lambda (g) \,dg\) for the left regular representation
    \(\Lambda\). This conditions can characterise the image \(\oper{W}_f
    V\) among all functions on \(G\) or \(X\).
  \item Often, the image of (induced) covariant transform has the
    reproducing property \( f(y)=\int_X f(x)\,\bar{k}_y(x)\,dx\), where
    \(k_y\) is the covariant transform of the shifted mother wavelet
    \(\uir{}{}(\map{s}(y)^{-1}) v\).
  \item For a representation \(\uir{}{}\) of \(G\) in a vector space
    \(V\) there is a representation \( \uir[3]{}{}(g_1,g_2): A
    \mapsto\uir{}{}(g_2)^{-1}A \uir{}{}(g_1)\) of \(G\times G\) on the
    space \(B(V)\) of bounded linear operators on \(V\).
  \item The covariant symbol \(
    \tilde{A}(g_1,g_2)=F(\uir[3]{}{}(g_1,g_2)
    A)=F(\uir{}{}(g_2)^{-1}A\uir{}{}(g_1))\) is the covariant transform
    defined by the representation \(\uir[3]{}{}\) and the operator \(F:
    B(V)\rightarrow U\).
  \end{itemize}
\item Contravariant transform.
  \begin{itemize}
  \item For a representation \(\uir{}{}\) of \(G\) and a summable
    function \(k\) on \(G\), the integrated representation is
    \(\uir{}{}(k)=\int_{G} k(g)\,\uir{}{}(g)\,dg\).
  \item The relative convolution is an integrated representation over a
    homogeneous space \( \uir{}{}(k)=\int_{X}
    k(x)\,\uir{}{}(\map{s}(x))\,dx\).
  \item The composition of two \(\uir{}{}(k_1) \uir{}{}(k_2)\) relative
    convolutions produces the twisted convolution \( \uir{}{}(k_1 \natural
    k_2)=\uir{}{}(k_1) \uir{}{}(k_2)\).
  \item For an integrated representation or relative convolution
    \(\uir{}{}\) and a fixed vector \(w\in V\) the contravariant transform
    of a function \(k\) is \( \oper{M}_w^{\uir{}{}} (k)=
    \uir{}{}(\check{k})w\), where \(\check{k}(g)=k(g^{-1})\).
  \item The twisted convolution of two wavelet transforms is
    \(\oper{W}_{v_1} u_1 \natural \oper{W}_{v_2} u_2= \theta(u_2,v_1)
    \oper{W}_{v_2} u_1\).
  \end{itemize}
\end{itemize}

\printindex

\small
\providecommand{\noopsort}[1]{} \providecommand{\printfirst}[2]{#1}
  \providecommand{\singleletter}[1]{#1} \providecommand{\switchargs}[2]{#2#1}
  \providecommand{\irm}{\textup{I}} \providecommand{\iirm}{\textup{II}}
  \providecommand{\vrm}{\textup{V}} \providecommand{\cprime}{'}
  \providecommand{\eprint}[2]{\texttt{#2}}
  \providecommand{\myeprint}[2]{\texttt{#2}}
  \providecommand{\arXiv}[1]{\myeprint{http://arXiv.org/abs/#1}{arXiv:#1}}
  \providecommand{\doi}[1]{\href{http://dx.doi.org/#1}{doi:
  #1}}\providecommand{\CPP}{\texttt{C++}}
  \providecommand{\NoWEB}{\texttt{noweb}}
  \providecommand{\MetaPost}{\texttt{Meta}\-\texttt{Post}}
  \providecommand{\GiNaC}{\textsf{GiNaC}}
  \providecommand{\pyGiNaC}{\textsf{pyGiNaC}}
  \providecommand{\Asymptote}{\texttt{Asymptote}}

\end{document}